\documentclass[reqno]{jgokova}
\usepackage{epsfig,graphicx}

\begin{document}

\newcommand{\ben}{\begin{enumerate}}
\newcommand{\een}{\end{enumerate}}
\newcommand{\be}{\begin{equation}}
\newcommand{\ee}{\end{equation}}
\newcommand{\bea}{\begin{eqnarray}}
\newcommand{\eea}{\end{eqnarray}}
\newcommand{\bc}{\begin{center}}
\newcommand{\ec}{\end{center}}

\newtheorem{thm}{Theorem}[section]
\newtheorem{cor}[thm]{Corollary}
\newtheorem{lem}[thm]{Lemma}
\newtheorem{prop}[thm]{Proposition}
\newtheorem{conj}[thm]{Conjecture}

\theoremstyle{definition}
\newtheorem{defn}[thm]{Definition}

\theoremstyle{remark}
\newtheorem{rem}[thm]{\rm\bfseries{Remark}}
\newtheorem*{notation}{Notation}

\newtheorem{ques}[thm]{\rm\bfseries{Question}}
\newtheorem{cons}[thm]{\rm\bfseries{Construction}}
\newtheorem{exm}[thm]{\rm\bfseries{Example}}


\setcounter{page}{1}
\volume{3}

\newcommand{\D}{{\mathbb D}}
\newcommand{\C}{{\mathbb C}}
\newcommand{\R}{{\mathbb R}}
\newcommand{\Z}{{\mathbb Z}}
\newcommand{\Q}{{\mathbb Q}}
\renewcommand{\P}{{\mathbb P}}
\newcommand{\s}{{\mathbb S}}
\newcommand{\B}{{\mathbb B}}
\newcommand{\I}{{\mathbb I}}
\newcommand{\h}{{\mathbb H}}
\newcommand{\e}{{\mathbb E}}

              \newcommand{\J}{{\mathcal J}}
              \newcommand{\M}{{\mathcal M}}
              \newcommand{\W}{{\mathcal W}}
              \newcommand{\U}{{\mathcal U}}
              \newcommand{\T}{{\mathcal T}}
              \newcommand{\V}{{\mathcal V}}
              \newcommand{\E}{{\mathcal E}}
              \newcommand{\F}{{\mathcal F}}
              \renewcommand{\L}{{\mathcal L}}
              \renewcommand{\O}{{\mathcal O}}
              \newcommand{\N}{{\mathcal N}}
              \newcommand{\G}{{\mathcal G}}
              \renewcommand{\H}{{\mathcal H}}
              \newcommand{\bb}{{\mathcal B}}

\title[Resolving symplectic orbifolds]{Resolving symplectic orbifolds with applications to finite group actions}
\author[CHEN]{Weimin Chen}

\thanks{}

\address{Department of Mathematics and Statistics, University of Massachusetts at Amherst, Amherst, MA 01003, USA}
\email{wchen@math.umass.edu}

\address{}
\email{}

\begin{abstract}
We associate to each symplectic $4$-orbifold $X$ a canonical smooth symplectic resolution
$\pi: \tilde{X}\rightarrow X$, which can be done equivariantly if $X$ comes with a symplectic
$G$-action by a finite group. Moreover, we show that the resolutions of the symplectic $4$-orbifolds 
$X/G$ and $\tilde{X}/G$ are in the same symplectic birational equivalence class; 
in fact, the resolution of $\tilde{X}/G$ can be reduced to that of $X/G$ by successively 
blowing down symplectic $(-1)$-spheres. 

To any finite symplectic $G$-action on a $4$-manifold $M$, we associate a pair $(M_G,D)$,
where $\pi: M_G\rightarrow M/G$ is the canonical resolution of the quotient orbifold and $D$ 
is the pre-image of the singular set of $M/G$ under $\pi$. We propose to study the group action 
on $M$ by analyzing the smooth or symplectic topology of $M_G$ as well as the embedding 
of $D$ in $M_G$. In this paper, an investigation on the symplectic Kodaira dimension $\kappa^s$ of
$M_G$ is initiated. In particular, we conjecture that $\kappa^s(M_G)\leq \kappa^s(M)$.  
The inequality is verified for several classes of symplectic $G$-actions, including any actions 
on a rational surface or a symplectic $4$-manifold with $\kappa^s=0$.
\end{abstract}
\keywords{Orbifold singularity, symplectic resolution, $4$-manifold, finite group action,
branched covering, configuration of symplectic surfaces, symplectic Kodaira dimension.}

\maketitle

\section{Introduction and the main results}
The purpose of this paper is twofold. On the one hand, we are concerned with the basic question
of resolving orbifold singularities in symplectic geometry. We show that any symplectic $4$-orbifold
admits a canonical symplectic resolution (a more precise description will follow) without 
imposing any conditions on the structure of its singular set. The second goal of this paper 
is to introduce some new ideas to the study of symplectic finite group actions on $4$-manifolds,
which are based on analyzing the symplectic resolution canonically associated to the symplectic quotient orbifold of the group action. In forthcoming papers we shall develop these ideas in the context of several natural problems in group actions and $4$-manifolds.

\subsection{Resolving symplectic orbifolds}
The problem of resolving symplectic singularities was posed by Gromov (cf. \cite{G}) in the 1980s. 
The first progress was made by McCarthy and Wolfson in \cite{McW1}, where the authors dealt with the case of isolated orbifold singularities of a $4$-dimensional space. They formulated a natural notion of symplectic resolution and gave a construction in this case. In a subsequent work the construction was extended to isolated algebraic singularities in a symplectic $4$-manifold (cf. \cite{McW2}). In both works the authors employed certain techniques of gluing symplectic manifolds along a certain type of hypersurfaces to construct the symplectic resolution. Their construction can be regarded as a symplectic analog of resolution of singularities in algebraic geometry.

Going beyond the case of isolated singularities in dimension $4$, Cavalcanti, Fern\'{a}ndez, and
Mu\~{n}oz \cite{CFM} constructed symplectic resolutions of isolated orbifold singularities in all dimensions. Their work made use of the algebraic resolution of orbifold singularities as 
McCarthy and Wolfson did in \cite{McW2}, but the construction of symplectic resolution was 
different. For non-isolated orbifold singularities, the construction of symplectic resolution requires
new ideas. In  \cite{NP1, NP2}, Niederkr\"{u}ger and Pasquotto constructed symplectic resolutions
of any symplectic orbifolds which arise in the symplectic reduction of a Hamiltonian torus action.
Their idea is to construct an auxiliary circle action in a neighborhood of the singularities with the
largest isotropy group, and subsequently use it to perform a symplectic cut of Lerman (cf. \cite{L}).
This procedure produces a symplectic orbifold of singularities with smaller isotropy groups, giving
rise to a ``partial resolution" of the original symplectic orbifold. A symplectic resolution is then 
obtained by a sequence of such partial resolutions. As an application, the authors gave a different method for constructing symplectic resolution of an isolated cyclic quotient singularity (cf. \cite{NP1}). 
It uses neither the algebraic resolution nor the gluing techniques as it was done in the works 
of  \cite{McW1}, \cite{McW2} or \cite{CFM}, but rather, it amounts to performing a sequence 
of weighted blow-ups (as described in \cite{Gd}). Constructing a symplectic resolution for an arbitrary 
symplectic orbifold remains an open problem (see \cite{MR} for some very recent development). 

In this paper, we give a construction of symplectic resolution for an arbitrary symplectic $4$-orbifold without imposing any conditions on the structure of its singular set. In order to explain the nature of our construction, we shall look at the analogous situation in the complex analytic category. To this end, we let $X$ be a $2$-dimensional complex orbifold (e.g., the quotient orbifold of a holomorphic finite group action on a complex surface). The singular set $\Sigma$ (i.e., the subset of points whose isotropy group is nontrivial) may contain points of complex co-dimension $1$ in $X$. However, these points are automatically non-singular in the underlying space of $X$. More precisely, if we denote by $|X|$ the underlying analytic space of the orbifold $X$, then the points in $\Sigma$ which have complex co-dimension $1$ are regular points of $|X|$. In fact, $|X|$ is a complex orbifold with at most isolated singularities. With this understood, one may take a resolution of the analytic space $|X|$ to serve as a resolution of the complex orbifold $X$ (such a resolution is unique if we require it to be a minimal resolution, cf. \cite{Lau}). Our construction of symplectic resolution amounts to carrying out an analogous consideration in the symplectic category. 

The following notations (concerning symplectic $4$-orbifolds) will be used throughout. 
Let $(X,\omega)$ be a symplectic $4$-orbifold, with singular set denoted by $\Sigma$, i.e.,
$$
\Sigma=\{p\in X| \mbox{the isotropy group $\Gamma_p$ is nontrivial}\}.
$$
If we fix an $\omega$-compatible (orbifold) almost complex structure $J$, and let $g_J$ be the 
corresponding Riemannian metric, then at each $p\in \Sigma$, the tangent space $T_p X$
can be identified with $\C^2$, with the action of $\Gamma_p$ on $T_p X$ given by a subgroup
of $U(2)$. Consequently, $\Sigma$ can be decomposed as a disjoint union 
$\Sigma^0\sqcup \Sigma^\ast \sqcup \Sigma^1$, where

\begin{itemize}
\item $\Sigma^0=\{p\in\Sigma| \mbox{the action of $\Gamma_p$ on $T_p X\setminus\{0\}$ is free}\}$.
\item $\Sigma^\ast=\{p\in\Sigma| \mbox{$\Gamma_p$ fixes a complex line in $(T_p X,J)$}\}$.
\item $\Sigma^1=\{p\in\Sigma| \mbox{the action of $\Gamma_p$ on $T_p X\setminus\{0\}$ is 
not free but is fixed-point free}\}$.
\end{itemize}

Both $\Sigma^0,\Sigma^1$ consist of finitely many points, but $\Sigma^\ast$ is a $2$-dimensional
smooth manifold such that $\omega|_{\Sigma^\ast}$ is an area form. We can compactify each connected component of $\Sigma^\ast$ in $X$ by adding points from $\Sigma^1$. Let $\{\Sigma_i\}$ be the set of compactified connected components of $\Sigma^\ast$. Then each $\Sigma_i$ is a symplectic orbifold surface in $(X,\omega)$ (possibly immersed), with the points of self-intersection 
of each $\Sigma_i$ and the points of intersection of distinct $\Sigma_i,\Sigma_j$ contained in 
$\Sigma^1$. We denote by $|X|$ the underlying space of $X$. 

In the symplectic case, it is not difficult to see that $|X|$ is a smooth orbifold with at most isolated
singularities (we shall explain this in more detail in Section 2). In particular, $|X|$ is non-singular
along the $2$-dimensional components $\Sigma^\ast$. However, the symplectic form $\omega$ 
is singular on $|X|$ along $\Sigma^\ast$. With this understood, the main step in the construction of the symplectic resolution of $(X,\omega)$ is to show that the orbifold $|X|$ supports a natural symplectic structure; in fact we shall de-singularize $\omega$ along $\Sigma^\ast$. (We remark 
that in the complex analytic situation, the complex structure is automatically non-singular along
$\Sigma^\ast$, so this step is not necessary.) 
For the purpose of applications in finite group actions, we shall give an equivariant version of this construction. The following theorem is the main technical result.

\begin{thm}
Let $(X,\omega)$ be a symplectic $4$-orbifold, and let $G$ be a finite group acting smoothly 
on the $4$-orbifold $X$, preserving the symplectic structure $\omega$.  There are $G$-invariant neighborhoods $U$ of $\Sigma^1$ in $|X|$, which can be taken arbitrarily small, such that for any choice of $U$, there is a $G$-invariant symplectic structure $\omega^\prime$ on the orbifold $|X|$, such that $\omega^\prime=\omega$ on $|X|\setminus (\Sigma^\ast\cup U)$
(as symplectic forms) and $\omega^\prime=\omega$ on $\Sigma^\ast\setminus U$ as area forms.
Each $\Sigma_i$ is a symplectic orbifold surface in $(|X|,\omega^\prime)$, which may be singular with respect to the smooth structure of the orbifold $|X|$. The self-intersections and singular points of each $\Sigma_i$ occur only at points in $\Sigma^1$, and there is a $G$-invariant, $\omega^\prime$-compatible, integrable almost complex structure on $U$ with respect to which each 
$\Sigma_i\cap U$ is a (genuine) holomorphic curve. 
\end{thm}

\begin{rem}
(1) It is interesting to compare the symplectic de-singularization in Theorem 1.1 with the usual holomorphic de-singularization. Let $p\in \Sigma^\ast\setminus U$ and let $m$ be the order of the 
isotropy group $\Gamma_p$. If we let $(\delta,\phi)$, $(\rho,\psi)$ be the natural 
polar coordinates on $X$ and $|X|$ in the normal direction at $p$, then the fact that the symplectic forms $\omega,\omega^\prime$ agree in the complement of $\Sigma^\ast\setminus U$ forces the polar coordinates to obey the following equations 
$$
\rho^2=\frac{1}{m}\cdot \delta^2, \; \psi=m\cdot \phi,
$$
as in the normal direction the symplectic forms are given by $\delta d\delta\wedge d\phi$ and 
$\rho d\rho\wedge d\psi$ respectively (see Remark 2.3 for more details). On the other hand, in the holomorphic de-singularization, if $z$ is the normal holomorphic coordinate on $X$, then $w=z^m$
is the normal holomorphic coordinate on $|X|$. In terms of normal polar coordinates, $(\delta,\phi)$, 
$(\rho,\psi)$ are related by the equations $\rho=\delta^m$, $\psi=m\cdot \phi$.

(2) Even though the symplectic structure $\omega^\prime$ on $|X|$ 
depends a priori on fixing a small neighborhood $U$ of $\Sigma^1$, the corresponding (orbifold) 
canonical line bundle $K_{\omega^\prime}$ is uniquely determined up to isomorphism. We
shall denote it by $K_{|X|}$.

(3) In principle, the procedure in Theorem 1.1 can be reversed, i.e., a symplectic structure on $|X|$ also determines a symplectic structure on $X$. We shall state and prove the following corollary to 
this effect in the context of finite group actions.
\end{rem}

\begin{cor}
Let $M$ be a smooth $4$-manifold and $B \subset M$ be a smoothly embedded surface 
(maybe disconnected). Suppose a smooth $4$-manifold $\hat{M}$ is constructed by taking 
a branched covering of $M$ along $B$, with $\Gamma$ being the group of deck transformations.
Let $G$ be a finite group acting smoothly on $\hat{M}$ extending the group of deck
transformations $\Gamma$ and inducing a smooth $G/\Gamma$-action on $M$.
Then the $G$-action on $\hat{M}$ preserves a symplectic structure 
if and only if $M$ admits a $G/\Gamma$-invariant symplectic structure with respect to which 
$B$ is symplectic. 
\end{cor}

The idea to construct non-standard group actions by taking branched coverings along a 
non-standard embedding of a co-dimension $2$ submanifold goes back to the work of Giffen 
\cite{Gif} which produced the first counterexamples to the generalized Smith conjecture. 
A much more elaborated version of this idea was used by Fintushel, Stern and Sunukjian 
\cite{FSS} to construct the first examples of topologically equivalent but smoothly non-equivalent 
finite cyclic actions on irreducible $4$-manifolds. In particular, they constructed exotic $\Z_2$, 
$\Z_3$, and $\Z_4$ actions on the $K3$ surface. Corollary 1.3 implies that none of these exotic 
actions can be made symplectic. For example, consider the exotic involutions which were 
constructed by taking a double branched covering of $\C\P^2$ along a smoothly embedded 
surface $B$, where $B$ is topologically isotopic but not smoothly isotopic to the sextic. 
If the exotic involutions were symplectic, then by Corollary 1.3 there is a symplectic structure 
on $\C\P^2$ with respect to which $B$ is symplectic. But a theorem of Shevchishin \cite{She}
(see also Siebert and Tian \cite{ST}) implies that $B$ is smoothly isotopic to the sextic, which is a contradiction. 
We remark that Corollary 1.3 is relevant to the following question which is currently open.

\begin{ques}
Does there exist a smooth finite group action on a K\"{a}hler surface which preserves a 
symplectic structure but not a complex structure (i.e., a symplectic finite group action not smoothly 
equivalent to a holomorphic action)?
\end{ques}

With Theorem 1.1 at hand, a symplectic resolution of the symplectic $4$-orbifold 
$(X,\omega)$ is obtained
by simply applying the results in \cite{McW1}, \cite{McW2}, \cite{CFM}, or \cite{NP1} to the symplectic
$4$-orbifold $(|X|,\omega^\prime)$. More specifically, we shall adopt the construction in \cite{CFM}.
The result is summarized in the following theorem.

\begin{thm}
Let $(X,\omega)$ be a symplectic $4$-orbifold. There exist a symplectic $4$-manifold, 
denoted by $\tilde{X}$ and called the {\em resolution} of $(X,\omega)$, and a 
continuous map $\pi: \tilde{X}\rightarrow X$, with the following significance.
\begin{itemize}
\item [{(1)}] The map $\pi: \tilde{X}\setminus \pi^{-1}(\Sigma)\rightarrow X\setminus 
\Sigma$ is a diffeomorphism, $\pi^{-1}(\Sigma^\ast)$ is a smoothly embedded surface in
$\tilde{X}$ such that $\pi: \pi^{-1}(\Sigma^\ast)\rightarrow \Sigma^\ast$ is a diffeomorphism, 
and for each $p\in \Sigma^0\sqcup\Sigma^1$ which is a singular point of $|X|$ , 
$\pi^{-1}(p)$ is a configuration of embedded two-spheres $\{S_i\}$, where $S_i,S_j$, $i\neq j$, are either disjoint or intersect transversely at a single point, 
and no three distinct $S_i$ intersect in one point. If $p$ is a smooth point of $|X|$, $\pi^{-1}(p)$
is a single point.

\item [{(2)}] There are neighborhoods $U$ of $\Sigma^0\sqcup\Sigma^1$, which can be taken
arbitrarily small, such that for any choice of $U$, there is a symplectic structure $\tilde{\omega}$ 
on $\tilde{X}$ such that 
$\pi: (\tilde{X}\setminus \pi^{-1}(U\cup\Sigma^\ast),\tilde{\omega})\rightarrow (X\setminus (U\cup \Sigma^\ast),\omega)$ 
is a symplectomorphism, $\pi^{-1}(\Sigma^\ast)$ is symplectic with $\tilde{\omega}=\pi^\ast \omega$ on $\pi^{-1}(\Sigma^\ast\setminus U)$ as area forms, 
and for each $p\in \Sigma^0\sqcup\Sigma^1$ which is a singular point in $|X|$, 
the components of $\pi^{-1}(p)$ are symplectic and intersect 
positively. Furthermore, there exists a neighborhood $U^\prime$ of 
$\Sigma^1$, where $U^\prime\subset U$, with an $\tilde{\omega}$-compatible, integrable almost 
complex structure on $U^\prime$ such that $U^\prime\cap \pi^{-1}(\Sigma)$ is given by (genuine)
holomorphic curves.
\item [{(3)}] Suppose a finite group $G$ acts on $X$ smoothly, preserving the symplectic structure 
$\omega$. Then there is a natural smooth $G$-action on $\tilde{X}$
such that $\pi:\tilde{X}\rightarrow X$ is $G$-equivariant,  and that $\tilde{\omega}$ 
can be made $G$-invariant if the neighborhood $U$
is chosen $G$-invariant. Furthermore, let $V,W$ be the resolutions of the symplectic 
$4$-orbifolds $(X/G,\omega)$ and $(\tilde{X}/G,\tilde{\omega})$ respectively. Then either
$W=V$, or $W$ can be reduced to $V$ by successively blowing down
symplectic $(-1)$-spheres. 
\end{itemize}
\end{thm}

\begin{rem}
(1) The diffeomorphism type of $\tilde{X}$ is uniquely determined by the smooth orbifold $X$; in fact, it is the smooth $4$-manifold obtained by replacing a neighborhood of each singular point of the smooth orbifold $|X|$ by the minimal resolution of the neighborhood. In particular, the symplectic 
Kodaira dimension of $\tilde{X}$ is well-defined and depends only on the smooth structure of $X$
(cf. \cite{Li}). 

(2) Even though the symplectic structure $\tilde{\omega}$ depends on a choice of a neighborhood of 
$\Sigma^0\sqcup \Sigma^1$, the canonical line bundle $K_{\tilde{\omega}}$ is uniquely determined up to isomorphism. 
We shall denote it by $K_{\tilde{X}}$. This said, the Gromov-Taubes invariant of 
$(\tilde{X},\tilde{\omega})$ is independent of the choice of $\tilde{\omega}$ (cf. \cite{T}).  
\end{rem}

\subsection{Finite group actions on $4$-manifolds}
Smooth finite group actions on $4$-manifolds remain poorly understood in general, despite the
tremendous advances in smooth $4$-manifold theory following the revolutionary work of 
Donaldson in the early 1980s (cf. \cite{D}). A key ingredient in understanding a finite group
action is the so-called fixed-point set structure. To be more precise, consider $M$ which is an 
oriented smooth $4$-manifold equipped with an orientation preserving smooth action
of a finite group $G$. For any $g\in G$, the fixed-point set $Fix(g)$ of $g$ consists of finitely
many isolated points and $2$-dimensional components. Crucial to the understanding of the
$G$-action is information about $Fix(g)$ and the induced representation of $g$ on the tangent
space $T_p M$ at a fixed point $p\in Fix(g)$, for any given $g\in G$ and $p\in Fix(g)$. This information is constrained 
by the various $G$-index theorems (e.g., Lefschetz fixed point theorem, $G$-signature theorem, etc.), 
however, for a general smooth $4$-manifold, the fixed-point set structure of a smooth finite 
group action remains poorly understood. 

Fixed-point set structure is intimately related to the various rigidity properties of the group action.
In this sense, locally linear topological actions are the most flexible ones. Indeed, 
Edmonds and Ewing \cite{EE} showed that for any pseudo-free, locally linear $\Z_p$-action 
of prime order on a simply connected $4$-manifold, the Lefschetz fixed point theorem and 
the $G$-signature theorem are almost the only constraints on the fixed-point set structure. 
(Pseudo-free, in this case, means that the fixed-point set of the action contains no $2$-dimensional components.) Therefore, any rigidity phenomenon of smooth actions (in comparison with locally linear actions) reflects existence of additional constraints on the fixed-point set structure of the group actions, and vice versa, any additional constraints on the fixed-point set structure may lead to certain rigidity properties of the group actions. 

Primary examples of smooth actions are provided by holomorphic actions on complex 
K\"{a}hler surfaces, which serve as a model and motivate the study of smooth actions. 
A well-known rigidity phenomenon of holomorphic actions is the so-called ``homological rigidity" of 
automorphisms of a $K3$ surface, i.e., a holomorphic automorphism must be trivial if 
the induced action on the $K3$ lattice is trivial (cf. \cite{BP}). 
Knowing that this fails to be true for locally linear actions on $K3$ surfaces (cf. \cite{E}), 
Edmonds asked whether there exist nontrivial smooth actions on $K3$ surfaces which are homologically trivial
(see Kirby's Problem List \cite{K}, Problem 4.124 (B)). Another example of rigidity properties of
holomorphic actions concerns the order of automorphism group of a minimal algebraic surface 
of general type (cf. \cite{Xiao}). 

A natural class of smooth actions which generalize the holomorphic actions on K\"{a}hler surfaces
is given by the symplectic finite group actions, and a central problem in this regard is to
what extent the rigidity properties of holomorphic actions can be extended to this class of
smooth actions. When $M$ is symplectic and the $G$-action preserves a symplectic structure, the pseudo-holomorphic curve theory and Taubes' seminal work on symplectic $4$-manifolds \cite{T} can be adapted to the equivariant setting which gave some powerful techniques to study the group action (see the survey articles \cite{ChenS1, ChenS2}, and the more recent papers \cite{C2, CLW}). In particular, these techniques revealed additional constraints on the fixed-point set structure of the group action, and allowed one to extend (partially) the aforementioned rigidity properties of holomorphic actions to symplectic finite group actions. For example, in joint work with Kwasik \cite{CK1}, the author showed that there are no symplectic finite group actions on 
the standard smooth $K3$ surface which act trivially on homology (it is an interesting open question as whether this continues to be true if the smooth structure of the $K3$ surface is exotic). For another example, the author investigated in \cite{C1} the problem of bounding the order of a symplectic 
finite group action, and partially extended the work of Xiao in \cite{Xiao} to the symplectic category.

With the preceding understood, we introduce in this paper some new construction to the study
of symplectic finite group actions. The idea is to associate to each symplectic finite group action a symplectic $4$-manifold which is canonically determined by the group action, i.e., the resolution of the quotient orbifold. (We shall call it the {\it resolution} of the group action or the {\it resolution} of the corresponding $G$-manifold.) 

To be more concrete, let $M$ be a symplectic $4$-manifold equipped with a finite symplectic 
$G$-action. Denote by $\pi: M_G\rightarrow M/G$ the resolution of the quotient orbifold $M/G$ 
from Theorem 1.5 and let $D:=\pi^{-1}(\Sigma)\subset M_G$ be the pre-image 
of the singular set of $M/G$, which is in general a configuration of symplectic surfaces in $M_G$.
We shall investigate the group action on $M$ through $M_G$ and the embedding of $D$ in $M_G$. Note that topologically, $M$ may be regarded as a branched covering of $M_G$ along $D$. To be more precise, let $\hat{M}$ be the $4$-manifold which is the regular $G$-branched covering of $M_G$ along $D$, then $\hat{M}$ and $M$ are related by a successive blowing-down of 
$(-1)$-spheres which is equivariant with respect to the $G$-actions on $\hat{M}$ and $M$. 
Particularly, this gives an unified point of view for the construction of symplectic finite group actions,
which has been lagging behind for progress. Finally, we remark that in some sense, $M_G$ can be regarded as a certain ``manifold approximation'' of the singular quotient orbifold $M/G$, and the subset $D\subset M_G$, which contains vital information about the fixed-point set structure of the $G$-action, is a substitute for the singular set of $M/G$. 

The idea of studying a finite automorphism group of a complex surface through the
resolution of the quotient space had appeared in the algebraic geometry literature, e.g., in 
Nikulin \cite{Nik}, Xiao \cite{Xiao, X3}, and the unpublished work of Xiao \cite{X2}. So our
new approach to the study of symplectic finite group actions on $4$-manifolds is an attempt to generalize this idea to the symplectic category. 

One basic invariant associated to a symplectic $4$-manifold is its symplectic Kodaira 
dimension (cf. \cite{Li}).
We may wonder (1) how the rigidity of a group action is seen through the symplectic Kodaira dimension
of $M_G$, and (2) how the symplectic Kodaira dimensions of $M_G$ and $M$ are related. 

For the second question, it is known that for a finite holomorphic $G$-action on a complex surface 
$M$, 
the Kodaira dimensions of $M_G$, $M$ obey the following inequality:
$$
\kappa(M_G)\leq \kappa(M).
$$
(This follows easily from the definition of Kodaira dimension.) As the symplectic Kodaira 
dimension coincides with the (complex) Kodaira dimension in the  K\"{a}hler case (cf. \cite{Li}), it
is natural to speculate that the corresponding inequality might be true in the symplectic category.

\begin{conj}
Let $\kappa^s$ stand for the symplectic Kodaira dimension. Then
$$
\kappa^s(M_G)\leq \kappa^s(M).
$$
\end{conj}

Conjecture 1.7 offers, in some very rough way, a measurement as how much symplectic 
finite group actions resemble holomorphic actions. 

\begin{rem}
(1) Conjecture 1.7 is known to be true in some special cases. For example, it is true for many
symplectic $G$-actions on $\C\P^2$ or Hirzebruch surfaces because the corresponding
group actions are smoothly equivalent to a holomorphic action (cf. \cite{C0, C00, ChenS2, C2}). 
On the other hand, if $M$ is a homotopy $K3$ surface with trivial 
canonical line bundle and the symplectic $G$-action satisfies $b^{+}_2(M/G)=3$, then 
$M_G$ is also a homotopy $K3$ surface with trivial canonical line bundle, cf. \cite{CK3}. 
(In this case, $M/G$ has only isolated singularities so the construction of $M_G$ is straightforward.) 
In this last example, $\kappa^s(M_G)=\kappa^s(M)=0$ so Conjecture 1.7 is true. 

(2) Theorem 1.5(3) allows an inductive approach toward Conjecture 1.7.
More precisely, suppose $H$ is a normal subgroup of $G$, with $K=G/H$ being the quotient group. 
Then by Theorem 1.5(3),  there is an induced symplectic $K$-action on the resolution $M_H$, 
and furthermore, $\kappa^s((M_H)_K)=\kappa^s(M_G)$ (cf. \cite{Li}). Hence, if one can show 
$\kappa^s(M_H)\leq \kappa^s(M)$ and $\kappa^s((M_H)_K)\leq \kappa^s(M_H)$, 
then one has $\kappa^s(M_G)\leq \kappa^s(M)$.
\end{rem}

Concerning the first question, the most interesting and important example of a connection between 
rigidity of a group action and the (symplectic) Kodaira dimension of the resolution is given in the work 
of Xiao \cite{Xiao} on the order of automorphism group of a minimal algebraic surface of general type. 
In that paper, Xiao showed that if $M$ is a minimal algebraic surface of general type and $G$ is its 
automorphism group, then the number $|G|/c_1(M)^2$ has an interesting correlation with the topology, 
particularly the Kodaira dimension of $M_G$: 
\begin{itemize}
\item [{(i)}] If $M_G$ is of general type, then $|G|/c_1(M)^2\leq 1$. 
\item [{(ii)}] If $\kappa(M_G)=1$, then $|G|/c_1(M)^2\leq 3$. 
\item [{(iii)}] If $M_G$ is a ruled surface over a curve of genus $\geq 2$, then $|G|/c_1(M)^2\leq 10.5$. 
\item [{(iv)}] If $\kappa(M_G)=0$, then $|G|/c_1(M)^2\leq 288$. 
\item [{(v)}] If $M_G$ is a rational surface or a ruled surface over an elliptic curve, 
$|G|/c_1(M)^2\leq c$
for some universal constant $c>0$ (the computation of $c$ is given separately in \cite{X1}).
\end{itemize}

In this paper, we give some further evidence for Conjecture 1.7. In particular, we have a pretty 
good understanding when $\kappa^s(M)=-\infty$ or $0$. (Recall that a symplectic $4$-manifold 
is rational or ruled if and only if $\kappa^s=-\infty$.)

\begin{thm}
Let $M$ be a symplectic $4$-manifold with a finite symplectic $G$-action. 
\begin{itemize}
\item [{(1)}] If $M$ is a rational surface, so is $M_G$. 
\item [{(2)}] If $\kappa^s(M)=0$, then either $\kappa^s(M_G)=0$, or $M_G$ is a rational surface, 
or a ruled surface over $T^2$. 
\end{itemize}
\end{thm}

We remark that, for the case of $\kappa^s(M)=0$, the verification of Conjecture 1.7 (in Theorem 1.9) does
not rely on any information a priori about the symplectic $G$-action on $M$, while for the case where 
$M$ is rational, it requires some nontrivial results concerning the equivariant symplectic cone of rational 
$G$-surfaces in \cite{CLW}. 

The case where $\kappa^s(M)=1$ is more intricate. We shall only 
consider in this paper the symplectic $G$-actions such that $b_2^{+}(M/G)>1$. Under this assumption, 
the equivariant version of the Seiberg-Witten-Taubes theory yields very strong constraints for 
the fixed-point set structure of
the action (cf. \cite{CK1, C1}). Note that $b_2^{+}(M_G)=b_2^{+}(M/G)>1$, so one always has 
$\kappa^s(M_G)\geq 0$ under this assumption. The results are summarized in the following
theorem. 

\begin{thm}
Let $M$ be a minimal symplectic $4$-manifold with $\kappa^s(M)=1$, equipped with a finite symplectic $G$-action, where $G=\Z_p$ is of prime order $p$. Furthermore, assume 
$b_2^{+}(M/G)>1$.  
\begin{itemize}
\item [{(1)}] $M_G$ has torsion canonical class if and only if the $2$-dimensional fixed components 
of $G$ consist of tori $\{T_i\}$ with self-intersection zero and the isolated fixed points of $G$ are 
all of type $(1,-1)$, and $c_1(K_M)=(p-1)\sum_i T_i$.
\item [{(2)}] In general, 
$c_1(K_{M_G})^2=-\frac{2(p-1)^2}{p}\cdot s+\sum_m K_m^2$, where $s$ is the
number of $(-2)$-spheres fixed by $G$, and $K_m$ denotes the canonical class of the minimal
resolution of the singular point of $M/G$ corresponding to the isolated fixed point $m$ of $G$.
\item [{(3)}] If the $G$-action is homologically trivial, then $c_1(K_{M_G})^2=-\sum_m \chi_m$,
where $\chi_m$ is the number of exceptional divisors in the minimal resolution of the singular point of $M/G$ corresponding to the isolated fixed point $m$ of $G$.
\item [{(4)}] If $c_1(K_{M_G})^2=0$, then $M_G$ must be minimal. 
\end{itemize}
\end{thm}

\begin{rem}
(1) When $M_G$ has torsion canonical class, the topology of $M_G$ is severely constrained 
(cf. \cite{Li, Li1,Bauer}). In particular, with $b_2^{+}(M_G)>1$, $M_G$ is either a $\Q$-homology $T^2$-bundle over $T^2$, or a homology $K3$ surface. In turn, this gives severe constraints
on the topology of $M$ and the order $p$ of the group $G$. To see this, let $F$ be the number
of isolated fixed points of $G$ (which in this case equals the Lefschetz number), then
$$
b_2^{-}(M/G)+(p-1)F=b_2^{-}(M_G).
$$
With this understood, observe that $b_2^{-}(M_G)=2,3$ or $19$. 

(2) The expression of $c_1(K_{M_G})^2$ given in part (2) implies that $c_1(K_{M_G})^2\leq 0$.
This is consistent with Conjecture 1.7, because if $c_1(K_{M_G})^2>0$ were true, then one 
would have $\kappa^s(M_G)=2$ which violates Conjecture 1.7. In fact, when $c_1(K_{M_G})^2=0$,
Conjecture 1.7 is true if and only if $M_G$ is minimal, which is confirmed in part (4). When
$c_1(K_{M_G})^2< 0$, $M_G$ must not be minimal by Taubes \cite{T}. In fact, $M_G$ should
contain at least $-c_1(K_{M_G})^2$ many $(-1)$-spheres, and Conjecture 1.7 is true if and only if $M_G$ contains exactly $-c_1(K_{M_G})^2$ many $(-1)$-spheres. 

(3) Observe that the expression of $c_1(K_{M_G})^2$ in part (2) consists of rational numbers,
so the identity gives rise to some integrability conditions on these rational numbers. However,
these integrability conditions do not yield any new constraints on the fixed-point set structure (see 
Remark 4.2 for more details).
\end{rem}

In general, $M_G$ contains a large number of $(-1)$-spheres, which must all intersect the 
subset $D$ if $M$ is assumed to be minimal. It is natural to attempt to extract new constraints
for the fixed-point set structure from the interaction of these  $(-1)$-spheres with $D$. We look
into this matter by examining the homologically trivial symplectic $\Z_3$-actions on a 
homotopy $K3$ surface. It turns out that in each case, we can see the $(-1)$-spheres in
$M_G$ explicitly, and blowing down these $(-1)$-spheres we obtain the minimal model of $M_G$.
For all these possible group actions, we have $\kappa^s(M_G)\leq 1$, so Conjecture 1.7 is true. 
On the other hand, our new approach by examining $M_G$ does not yield any new constraints. 
See Example 4.3 for more details. 


The organization of this paper is as follows. Section 2 is devoted to the proof of the main
technical result, Theorem 1.1. After a preliminary lemma on smooth finite group actions on a 
general smooth orbifold, the section proceeds to a detailed discussion on how to 
de-singularize the symplectic structure along the $2$-dimensional singular strata of the symplectic $4$-orbifold in question, which forms the bulk of the technical work of this paper. Section 2 ends 
with the proof of Corollary 1.3. The construction of symplectic resolution (Proof of Theorem 1.5) 
is given in Section 3, which also contains two propositions (Propositions 3.2 and 3.3): one concerning the canonical class of the symplectic resolution, and the other concerning equivariant blowing down. The final section, Section 4, is mainly devoted to the proofs of Theorems 1.9 and 1.10 (on the symplectic Kodaira dimension of $M_G$), however, we also include 
at the end of the section a discussion of homologically trivial symplectic $\Z_3$-actions on 
a homotopy $K3$ surface from the point of view of symplectic resolution. 

\section{Symplectic orbifold structure on the underlying space}
This section is devoted to the proof of Theorem 1.1. First of all, we shall describe a natural smooth
orbifold structure on $|X|$, which has at most isolated singularities. With respect to this smooth orbifold structure, each smooth point in $X$, i.e., each point in $X\setminus \Sigma$, is a smooth
point in $|X|$, and each point of $\Sigma^0$ is a singular point of $|X|$ with the same orbifold local 
chart. As for each point $p\in \Sigma^\ast$, there is a smooth local chart of complex coordinates 
$(w_1,w_2)$ of $p$ in $|X|$ with the following property: there is a smooth local orbifold chart of 
$p$ in $X$ given by complex coordinates $(z_1,z_2)$, with a smooth $\Z_m$-action generated 
by $(z_1,z_2)\mapsto (z_1,\exp (2\pi i/m) z_2)$, such that $w_1=z_1$ and $w_2=z_2^m$. In particular, each point $p\in \Sigma^\ast$ is a smooth point in $|X|$ and $\Sigma^\ast$ is smoothly 
embedded in $|X|$. Finally, for each point 
$p\in \Sigma^1$, a local orbifold chart of $p$ in $|X|$ is obtained as follows. We take an orbifold 
chart $(\R^4,\Gamma_p)$ of $p$ in $X$, where $\Gamma_p$ acts on $\R^4$ linearly. 
Let $\s^3$ be the unit sphere in $\R^4$. Then $\s^3/\Gamma_p$ is a $3$-orbifold whose singular 
set is a link denoted by $L$. The underlying space $|\s^3/\Gamma_p|$ is a $3$-manifold, 
and a neighborhood of $p$ in $|X|$ is given by a cone over $|\s^3/\Gamma_p|$. With this understood, let $N_p$ be the normal subgroup of $\Gamma_p$ generated by the isotropy 
subgroups of the components of the singular set $L$ of the $3$-orbifold $\s^3/\Gamma_p$, 
and let $\Gamma_p^\prime:=\Gamma_p/N_p$. Then the fundamental group of the $3$-manifold 
$|\s^3/\Gamma_p|$ is isomorphic to $\Gamma_p^\prime$, which is obviously a finite group. The universal cover of $|\s^3/\Gamma_p|$ is $\s^3$, and $|\s^3/\Gamma_p|$ is diffeomorphic to the quotient of $\s^3$ by a linear action of $\Gamma_p^\prime$. When $\Gamma_p^\prime$ is nontrivial, the point $p\in |X|$ has a natural orbifold chart given by $(\R^4,\Gamma_p^\prime)$.
If $\Gamma_p^\prime$ is trivial, then $p\in |X|$ is a smooth point. 

\subsection{Finite group actions on symplectic orbifolds}

We start off with a preliminary lemma concerning symplectic finite group actions on a symplectic orbifold. We include a detailed proof here for completeness; this type of arguments will also be used
in the proof of Lemma 2.2 and we shall be brief there. For a general reference on smooth orbifolds 
we refer the reader to \cite{C01}.

\begin{lem}
Let $X$ be a smooth $n$-orbifold which is equipped with a smooth finite group action by $G$. Then
the quotient space $X/G$ is naturally a smooth $n$-orbifold. Furthermore, if $\omega$ is a symplectic structure on $X$
which is preserved under the $G$-action, then $\omega$ descends to a symplectic structure on the orbifold $X/G$,
making it naturally a symplectic orbifold.
\end{lem}

\begin{proof}
For any $p\in X$, let $G_p$ be the isotropy subgroup at $p$, i.e., 
$$
G_p:=\{g\in G| g\cdot p=p\},
$$
and denote by $\bar{p}$ the image of $p$ in the quotient space. Furthermore, 
let $(\R^n,\Gamma_p)$ be a local orbifold chart of $X$ at $p$. With this understood, if $G_p$ is trivial, then $G\cdot p$ is
a free orbit, so that a neighborhood of $p$ is mapped homeomorphically onto a neighborhood of $\bar{p}$ in the quotient space.
Hence in this case, there is a natural smooth orbifold chart of the quotient space at $\bar{p}$, i.e., $(\R^n,\Gamma_p)$, which is
independent of the choice of $p$ up to equivalence. 

Suppose $G_p$ is nontrivial. Then there is a $\Gamma_p$-invariant open subset $U_p$ of $\R^n$ containing $0$, 
such that $U_p/\Gamma_p$ is a $G_p$-invariant neighborhood of $p$. Furthermore, for any $g\in G_p$, there is a lifting
of $g$ to a diffeomorphism $\tilde{g}: U_p\rightarrow U_p$, and any two such liftings $\tilde{g},\tilde{g}^\prime$ differ by an
element of $\Gamma_p$, i.e., $\tilde{g}^\prime=h\circ \tilde{g}$ for some $h\in \Gamma_p$. With this understood, for each
$g\in G_p$ we fix a lifting $\tilde{g}: U_p\rightarrow U_p$. Note that for any $g_1,g_2\in G_p$, the lifting of $g_1g_2$ differs 
from $\tilde{g}_1\circ \tilde{g}_2$ by an element of $\Gamma_p$. Now we let $\tilde{\Gamma}_p$ be the group of 
self-diffeomorphisms of $U_p$ generated by $h$ and $\tilde{g}$ for all $h\in \Gamma_p$ and $g\in G_p$. Then the above 
property implies that the map from $\tilde{\Gamma}_p$ to $G_p$ sending each $h\in \Gamma_p$ to $1\in G_p$ and each
$\tilde{g}$ to $g$ is a surjective homomorphism, whose kernel is $\Gamma_p$. In particular, 
$\tilde{\Gamma}_p$ is a finite group, with a natural short exact sequence 
$$
1\rightarrow \Gamma_p\rightarrow \tilde{\Gamma}_p\rightarrow G_p\rightarrow 1.
$$
Observing that $U_p/\tilde{\Gamma}_p=(U_p/\Gamma_p)/G_p$, we can define $(U_p, \tilde{\Gamma}_p)$ to be a smooth 
orbifold chart of the quotient $X/G$ at $\bar{p}$. It is easy to see that it is independent of the choice of $p$ up to equivalence. 
This puts a natural smooth orbifold structure on $X/G$. 

Let $\eta$ be a differential form on $X$. For any $g\in G$, the pullback $g^\ast\eta$ is defined as follows. Let $p,q\in X$
such that $g\cdot p=q$. Let $(U_p,\Gamma_p)$, $(U_q,\Gamma_q)$ be local orbifold charts at $p,q$ respectively,
such that a lifting of $g$, i.e., a diffeomorphism $\tilde{g}:U_p\rightarrow U_q$ covering $g: U_p/\Gamma_p\rightarrow
U_q/\Gamma_q$, exists. Furthermore, let $\eta_q$ be the local $\Gamma_q$-invariant differential form on $U_q$ 
representing $\eta$. Then the pullback $\tilde{g}^\ast \eta_q$ is a $\Gamma_p$-invariant differential form on $U_p$, which
is independent of the choice of the local lifting $\tilde{g}$ of $g$. This is because for any lifting $\tilde{g}^\prime: U_p\rightarrow U_q$,
there is a $h\in \Gamma_q$ such that $\tilde{g}^\prime=h\circ \tilde{g}$, and in particular, for any $k\in \Gamma_p$, 
$\tilde{g}\circ k=k^\prime\circ \tilde{g}$ for some $k^\prime\in \Gamma_q$. With this understood, the pullback $g^\ast\eta$ is
the differential form on $X$ determined by the local $\Gamma_p$-invariant forms $\tilde{g}^\ast \eta_q$. The above description
of pullback forms immediately implies that if $\omega$ is a symplectic form on $X$ preserved by the $G$-action, then for any
$p\in X$ and $g\in G_p$, $\tilde{g}^\ast\omega_p=\omega_p$, where $\omega_p$ is the $\Gamma_p$-invariant symplectic
form on $U_p$ representing $\omega$. It follows easily that each $\omega_p$ is also $\tilde{\Gamma}_p$-invariant, and 
hence $\omega$ descends to a symplectic form on the quotient orbifold $X/G$. This finishes the proof of the lemma.

\end{proof}

We remark that in the case of $4$-orbifolds, if $G$ acts on $X$ smoothly, then the induced 
$G$-action on $|X|$ is a smooth action of orbifolds, and furthermore, $|X/G|=||X|/G|$ as smooth orbifolds.

\subsection{De-singularizing symplectic forms along co-dimension $2$ singular strata}
The key issue in  the proof of Theorem 1.1 is to de-singularize the symplectic structure $\omega$ of $X$ along the subset $\Sigma^\ast$ of the singular set, as $\Sigma^\ast$ lies in the smooth locus of $|X|$. Recall that $\Sigma_i$ denotes the closure of a connected component of $\Sigma^\ast$ in
$\Sigma$. We shall consider first the special case where $\Sigma_i$ lies entirely in $\Sigma^\ast$,
and then explain how to extend the argument to the general case. To this end, let $\Sigma_i$
be such a component. To ease the notation, we denote $\Sigma_i$ by $S$, which is a compact 
closed Riemann surface. Let $H$ be the subgroup of $G$ which leaves $S$ invariant.

\begin{lem}
The symplectic structure $\omega$ can be de-singularized along $S$, i.e., there exists a symplectic structure $\omega^\prime$ of $|X|$ defined in a neighborhood of $S$, such that (1) $\omega^\prime$
is $H$-invariant, (2) $\omega^\prime=\omega$ in the complement of $S$, and (3) 
$\omega^\prime|_S=\omega|_S$ as area forms. 
\end{lem}

\begin{proof}
We begin by noting that for each $p\in\Sigma_i$, the isotropy group $\Gamma_p$ is cyclic of order $m>1$ where $m$ is independent of $p$. Let $\nu\rightarrow S$ be the normal bundle of $S$ in the 
orbifold $X$, which is defined to be the quotient bundle of $TX|_S$ by the sub-bundle $TS$, 
and as such, it comes with a structure of a $\Z_m$-complex line bundle over $S$ once we fix 
an $\omega$-compatible (orbifold) almost complex structure $J$. The $\Z_m$-action on the fibers 
of $\nu$ is given by the complex multiplication by an $m$-th root of unity, and the $\Z_m$-action on the base $S$ is trivial. 
Furthermore, the corresponding Riemannian metric $g_J$ gives rise to a metric on $\nu$. 
Let $D(\nu,r)$ be the disc-bundle of $\nu$ of radius $r$. Then via the exponential map associated to 
$g_J$, $(D(\nu,r),\Z_m)$ gives a global orbifold chart of $X$ near $S$ for small $r>0$ ; in particular, 
a neighborhood of $S$ in the underlying space $|X|$ is given by $D(\nu,r)/\Z_m$ homeomorphically.

\vspace{2mm}

{\bf Claim:} {\it There is a finite group $\tilde{H}$, with a short exact sequence 
$$
1\rightarrow \Z_m\rightarrow \tilde{H}\rightarrow H\rightarrow 1,
$$
such that the $H$-action near $S$ lifts to a $\tilde{H}$-action on the global orbifold chart 
near $S$. Furthermore, $\omega$ is $\tilde{H}$-invariant. 
}

\vspace{2mm}

{\bf Proof of Claim:} By averaging the metric $g_J$, we may assume it is $G$-invariant. With this assumption, the neighborhood $D(\nu,r)/\Z_m$ of $S$ is $H$-invariant; in particular, the boundary 
$\partial (D(\nu,r)/\Z_m)=\partial D(\nu,r)/\Z_m$ is invariant under the $H$-action. Now note that 
$\partial D(\nu,r)\rightarrow \partial D(\nu,r)/\Z_m$ is an $m$-fold regular cyclic covering map. It follows easily that for each $h\in H$, $h:\partial D(\nu,r)/\Z_m\rightarrow \partial D(\nu,r)/\Z_m$
can be lifted to a diffeomorphism $\tilde{h}:\partial D(\nu,r)\rightarrow \partial D(\nu,r)$. Then a similar argument as in the proof of Lemma 2.1 shows that the elements $\xi\in\Z_m$ and $\tilde{h}$, 
$h\in H$, generate a finite group $\tilde{H}$, with a natural short exact sequence 
$1\rightarrow \Z_m\rightarrow \tilde{H}\rightarrow H\rightarrow 1$. Furthermore, there is a smooth
action of $\tilde{H}$ on $D(\nu,r)$ inducing the $H$-action on the neighborhood $D(\nu,r)/\Z_m$ of $S$, and $\omega$ is preserved by the $\tilde{H}$-action. Hence the claim.

\vspace{2mm}

We continue with the proof of Lemma 2.2. Since $g_J$ is $G$-invariant, so is $J$, so that 
the normal bundle $\nu$ is a complex $\tilde{H}$-line bundle. Let $\pi: Y\rightarrow S$ be the principal $\s^1$-bundle associated to the complex line bundle $\nu$, which comes with a natural smooth bundle action by $\tilde{H}$. Note that the $\Z_m$-action on $Y$ induced by $\Z_m\rightarrow \tilde{H}$ is simply given by the inclusion of $\Z_m$ into $\s^1$. 
Let $Y^\prime:=Y/\Z_m$ be the quotient, and let $\pi^\prime: Y^\prime\rightarrow S$ 
be the corresponding principal $\s^1$-bundle. Then there is a natural smooth $H$-action 
on the principal $\s^1$-bundle $Y^\prime$. With this understood, note that $D(\nu,r)/\Z_m$ is 
$H$-equivariantly homeomorphic to a disc-bundle $D(\nu^\prime)$ of the complex line bundle 
$\nu^\prime\rightarrow S$ associated to the principal $\s^1$-bundle $Y^\prime$. Furthermore,
with respect to the (orbifold) smooth structure of $|X|$, a neighborhood of $S$ is given by 
$D(\nu^\prime)$ diffeomorphically.

We shall regard $D(\nu^\prime)$ as a smooth chart on $|X|$ near $S$, and as such,
$\omega$ may be regarded as a symplectic structure in the complement 
of the zero section of $D(\nu^\prime)$. This said, we shall describe next how to construct a 
symplectic structure $\omega^\prime$ on $D(\nu^\prime)$, which is $H$-invariant and agrees 
with $\omega$ in the complement of the zero section of $D(\nu^\prime)$. 

To this end, note 
first that $\omega$ defines an $\tilde{H}$-invariant symplectic structure on $D(\nu,r)$ via pull-back.
We shall begin by describing a model for $\omega$ on $D(\nu,r)$. Recall that $\pi:Y\rightarrow S$ 
is the $\tilde{H}$-equivariant principal $\s^1$-bundle associated to $\nu$. We pick an $\tilde{H}$-invariant connection $1$-form $\alpha$ on $Y$ and let $\kappa$ be the $2$-form on $S$ such that 
$\pi^\ast \kappa=-d\alpha$. Finally, let $\eta$ be the area form on $S$ which is the restriction of 
$\omega$ on $S$. Then for sufficiently small $r>0$, the following $2$-form is symplectic on 
$Y\times (-r,r)$:
$$
\omega_0:=\pi^\ast(\eta+ t\kappa)+\alpha\wedge dt,
$$
where $t$ is the coordinate on $(-r,r)$. The symplectic form $\omega_0$ is clearly $\tilde{H}$-invariant, 
where the $\tilde{H}$-action on $Y\times (-r,r)$ is trivial in the last factor. 

To relate the symplectic structures on $Y\times (-r,r)$ and $D(\nu,r)$, we employ a technique called 
symplectic cutting due to Lerman (cf. \cite{L}). To this end, note that the natural $\s^1$-action on 
$Y\times (-r,r)$ is Hamiltonian with respect to the symplectic structure $\omega_0$, with a 
Hamiltonian function given by $h: (y,t)\in Y\times (-r,r)\mapsto t$. In particular, each $Y\times \{t\}$ 
is a level set of $h$. To describe the procedure of symplectic cutting, consider the symplectic 
manifold $(Y\times (-r,r))\times \C$ with the symplectic structure 
$\omega_0 \oplus \frac{i}{2}dz\wedge d\bar{z}$ and the Hamiltonian $\s^1$-action given by 
$$
\lambda \cdot ((y,t),z)=((\lambda\cdot y,t),\lambda\cdot z), \;\; \forall \lambda\in\s^1.
$$
A Hamiltonian function of the $\s^1$-action is given by $\H((y,t),z)=h(y,t)-\frac{1}{2}|z|^2$. 
One can easily verify that $0$ is a regular value of $\H$. The level set $\H^{-1}(0)$ is easily 
seen to be the subset 
$$
\bigsqcup_{0\leq t<r} Y\times \{t\}\times \{|z|^2=2t\}.
$$
With this understood, the symplectic reduction $\H^{-1}(0)/\s^1$ can be easily identified
with $D(\nu,r)$, under which $Y\times (0,r)\subset \H^{-1}(0)/\s^1$ is identified with the 
complement of the zero section in $D(\nu,r)$. Furthermore, the symplectic structure on 
$\H^{-1}(0)/\s^1$ agrees with the symplectic structure $\omega_0$ on $Y\times (0,r)$. 
For simplicity, we will continue to denote the symplectic structure on $D(\nu,r)=\H^{-1}(0)/\s^1$ 
by $\omega_0$. Clearly, this symplectic form on $D(\nu,r)$ is $\tilde{H}$-invariant.

Note that on the zero section $S\subset D(\nu,r)$, $\omega_0=\eta=\omega$. Hence by
the equivariant version of the Weinstein neighborhood theorem, $\omega$ is 
$\tilde{H}$-equivariantly symplectomorphic to $\omega_0$ after taking $r>0$ sufficiently small. 
In other words, $\omega_0$ can serve as a model for $\omega$ in a neighborhood of $S$ in the
orbifold $X$.

With the preceding understood, denote by $D^\ast(\nu,r), D^\ast(\nu^\prime)$ the complement 
of the zero section in $D(\nu,r), D(\nu^\prime)$ respectively. The $\Z_m$-action on $D^\ast(\nu,r)$
is free, so that under the $H$-equivariant homeomorphism between $D(\nu,r)/\Z_m$ and
$D(\nu^\prime)$, $D^\ast(\nu,r)/\Z_m$ is sent to $D^\ast(\nu^\prime)$ $H$-equivariantly 
by a diffeomorphism. This said, the symplectic form $\omega_0$ on $D^\ast(\nu,r)$, which is
$\tilde{H}$-invariant, descends to an $H$-invariant symplectic form on $D^\ast(\nu,r)/\Z_m$. 
Denote the corresponding symplectic form on $D^\ast(\nu^\prime)$ by $\omega_0^\prime$. 

In order to extend $\omega_0^\prime$ to the disc-bundle $D(\nu^\prime)$, we shall again
employ the technique of symplectic cutting. To this end, we identify $D^\ast(\nu^\prime)$
with $Y^\prime\times (0,r)$, where $Y^\prime=Y/\Z_m$ is the principal $\s^1$-bundle 
associated to $\nu^\prime$, and consider $\omega_0^\prime$ as a symplectic form on
$Y^\prime\times (0,r)$. To identify $\omega_0^\prime$, note that the connection $1$-form
$\alpha$ on $Y$ descends to a $1$-form on $Y^\prime$. Furthermore, there is a connection 
$1$-form $\alpha^\prime$ on $Y^\prime$ such that $\alpha=\frac{1}{m}\alpha^\prime$.
Let $\kappa^\prime$ be the $2$-form on $S$ such that 
$(\pi^\prime)^\ast \kappa^\prime=-d\alpha^\prime$, where $\pi^\prime:Y^\prime\rightarrow S$.
Then $\kappa=\frac{1}{m}\kappa^\prime$. With this understood, the symplectic form 
$\omega_0^\prime$, as the descendent of $\omega_0$ under the free $\Z_m$-action, can be 
written on $Y^\prime\times (0,r)$ as
$$
\omega_0^\prime=(\pi^\prime)^\ast (\eta+\frac{t}{m}\kappa^\prime)
+\alpha^\prime\wedge d(\frac{t}{m}),
$$
where $t$ is the coordinate on $(0,r)$. This said, note that 
$\omega_0^\prime=(\pi^\prime)^\ast (\eta+\frac{t}{m}\kappa^\prime)
+\alpha^\prime\wedge d(\frac{t}{m})$ is in fact defined on the entire $Y^\prime\times (-r,r)$ and
is a symplectic form on it as long as $r>0$ is sufficiently small. 
Furthermore, note that the $\s^1$-action on $Y^\prime\times (-r,r)$
is Hamiltonian with respect to $\omega_0^\prime$, with a Hamiltonian function given by
$h^\prime: (y^\prime,t)\mapsto \frac{t}{m}$. A symplectic cutting procedure as we described
earlier, done $H$-equivariantly at the regular value $0$ of $h^\prime$, gives rise to an 
$H$-invariant symplectic form $\omega^\prime$ on $D(\nu^\prime)$, extending 
the symplectic form $\omega_0^\prime$ on $D^\ast(\nu^\prime)$. Note that $\omega^\prime=
\eta=\omega$ on the zero section as area forms. Now if we identify a neighborhood of $S$
in $|X|$ with $D(\nu^\prime)$, we obtained the symplectic structure $\omega^\prime$ with desired properties. This finishes the proof.

\end{proof}

\begin{rem}
It is somewhat a disturbing statement that the symplectic forms $\omega,\omega^\prime$ agree
on an open dense subset. We shall examine this more closely in local coordinates. Let 
$(\delta,\phi)$, $(\rho,\psi)$ be the polar coordinates in the fiber direction on $D(\nu,r)$ and 
$D(\nu^\prime)$, such that the restrictions of $\omega_0,\omega_0^\prime$ on the fibers are given by $\delta d\delta\wedge d\phi$ and $\rho d\rho\wedge d\psi$ respectively. Then if we follow 
through the symplectic cutting procedures in the proof of Lemma 2.2, it is easy to see that the polar coordinates are related by the equations 
$$
\rho^2=\frac{1}{m}\cdot \delta^2, \; \psi=m\cdot \phi,
$$
where $m$ is the order of the isotropy groups along $S$. Under the above relations, the forms 
$\delta d\delta\wedge d\phi$ and $\rho d\rho\wedge d\psi$ agree in the complement of the zero sections.

\end{rem}

\vspace{3mm}

\noindent{\bf Proof of Theorem 1.1}

\vspace{3mm}

If $\Sigma^1$ is empty, then Theorem 1.1 follows immediately from Lemma 2.2 with $U=
\emptyset$. Suppose $\Sigma^1\neq \emptyset$, and let $\Sigma_i$ be a component which
contains points in $\Sigma^1$. We shall extend the argument of Lemma 2.2 to first
de-singularize the symplectic structure $\omega$ along $\Sigma_i\setminus U$, where 
$U$ is a certain $G$-invariant neighborhood of $\Sigma^1$ which can be chosen arbitrarily small.
We summarize the result in the following lemma.

\begin{lem}
There are $G$-invariant neighborhoods $U$ of $\Sigma^1$ in $|X|$, which can be taken 
arbitrarily small, such that for any choice of $U$, there is a $G$-invariant symplectic structure
$\omega^\prime$ on $|X|\setminus U$, with $\omega^\prime=\omega$ in the complement of
$\Sigma^\ast\setminus U$ (as symplectic forms) and $\omega^\prime=\omega$ on 
$\Sigma^\ast\setminus U$ as area forms. 
\end{lem}

\begin{proof}
Note that $\Sigma_i\setminus U$ is a compact Riemann surface with boundary. Crucial to
extending the argument of Lemma 2.2 to the general situation is to describe an appropriate 
model for the symplectic structure $\omega$ near the boundary components of 
$\Sigma_i\setminus U$. We shall first explain how to choose the neighborhood $U$.

We begin by fixing a standard model for the symplectic structure $\omega$ near each point $p\in\Sigma^1$. 
Let $G_p$ be the subgroup of $G$ fixing $p$, $\tilde{\Gamma}_p$ be the isotropy group of the orbifold $X/G$ at $p$, 
and $\Gamma_p$ be the isotropy group of $X$ at $p$. Then these three groups fit into a short exact sequence 
(see the proof of Lemma 2.1)
$$
1\rightarrow \Gamma_p\rightarrow \tilde{\Gamma}_p\rightarrow G_p\rightarrow 1.
$$
By the equivariant Darboux theorem, there is an action of $\tilde{\Gamma}_p$ on $\C^2$
as a subgroup of $U(2)$, preserving the standard symplectic structure $\omega_0$ on
$\C^2$, such that for some $\epsilon_0>0$ which is independent of $p$, the induced action of 
$\Gamma_p$ on $(B^4(\epsilon_0),\omega_0)$ provides a model for an orbifold chart of $X$
at $p$, and the induced action of $G_p$ on the quotient $(B^4(\epsilon_0)/\Gamma_p,\omega_0)$ provides a model for the action of $G_p$ in the neighborhood of $p$ in $X$ (here $B^4(r)\subset \C^2$ is the open ball of radius $r$ with respect to the standard metric).

With the preceding understood, the $G$-invariant neighborhood $U$ will be taken to
be the union $\bigsqcup_{p\in \Sigma^1}U_p$, where $U_p$ is a neighborhood of $p$ 
modeled by $(B^4(t_0)/\Gamma_p,\omega_0)$ for some choice of $t_0<\epsilon_0$.
Fixing a $p\in \Sigma^1$, we let $\s^3(t)\subset B^4(\epsilon_0)$ be the sphere of radius $t$.
Since the action of $\Gamma_p$ is by a subgroup of $U(2)$, the Hopf fibration on $\s^3(t)$
is preserved under the action, and the Hopf fibration descends to a Seifert fibration on the 
$3$-orbifold $\s^3(t)/\Gamma_p$. Note that the singular set of $\s^3(t)/\Gamma_p$ is a union
of singular fibers of the Seifert fibration. Finally, for any $\Sigma_i$ which is not entirely
contained in $\Sigma^\ast$, if $p\in \Sigma_i$, then $\Sigma_i$ intersects with 
$\s^3(t)/\Gamma_p$ at a union of singular components. 

Now we fix a $t_0$ and let $\gamma$ be a singular component of the $3$-orbifold 
$\s^3(t_0)/\Gamma_p$. We shall describe a model for the symplectic structure $\omega$ 
near $\gamma$. Without loss of generality, we may assume that the pre-image of $\gamma$ 
in $\s^3(t_0)$ is given by $z_1=0$. A neighborhood of $\{z_1=0\}\cap \s^3(t_0)$ in 
$B^4(\epsilon_0)$, denoted by $W$, can be parametrized by coordinates $(r,\theta,\phi,t)$, 
where $0\leq r<r_0$, $\theta, \phi\in\R/2\pi\Z$, and $t$ lies in a small interval containing $t_0$, 
by the following map 
$$
(z_1,z_2)=(\frac{rt}{\sqrt{1+r^2}} e^{i(\theta+\phi)}, \frac{t}{\sqrt{1+r^2}} e^{i\phi}).
$$
With this understood, note that $\omega_0=r_1dr_1\wedge d\theta_1+r_2dr_2\wedge d\theta_2$ 
in polar coordinates (here $z_1=r_1e^{i\theta_1}$, $z_2=r_2e^{i\theta_2}$), hence in the new coordinates, a simple calculation shows that
$$
\omega_0=\frac{r^2t}{1+r^2} dt\wedge d\theta+\frac{t^2r}{(1+r^2)^2}dr\wedge d\theta
+tdt\wedge d\phi.
$$
We point out that $(r,\theta,\phi)$ gives a trivialization of the Hopf fibration near $z_1=0$ in
$\s^3(t_0)$, with $(r,\theta)$ for the base and $\phi$ for the fiber.

Let $\Gamma$ be the subgroup of $\Gamma_p$ which leaves the fiber $\{z_1=0\}\cap
\s^3(t_0)$ invariant. Then $\Gamma$ must be a cyclic group which acts on the neighborhood 
$W$ of $\{z_1=0\}\cap\s^3(t_0)$ by translations in the $\theta$ and $\phi$ variables. 
Furthermore, let $\Gamma^\prime$ be the subgroup of $\Gamma$ which leaves each fiber of the
Hopf fibration invariant, i.e., $\Gamma^\prime$ is the subgroup which acts by translation in the
$\phi$ variable only. We take the quotient by $\Gamma^\prime$ first, and by making a change of
variable $\phi\mapsto \frac{1}{|\Gamma^\prime|}\phi$ so that we continue to have $\phi\in\R/2\pi\Z$, we have the following expression for $\omega_0$ on the quotient $W/\Gamma^\prime$:
$$
\omega_0=\frac{r^2t}{1+r^2} dt\wedge d\theta+\frac{t^2r}{(1+r^2)^2}dr\wedge d\theta+
\frac{1}{|\Gamma^\prime|}tdt\wedge d\phi.
$$
With this understood, let $a$ be the order of the quotient group $\Gamma/\Gamma^\prime$. Then there is a generator $g\in\Gamma/\Gamma^\prime$ whose action on $W/\Gamma^\prime$ is given 
by
$$
g\cdot (r,\theta,\phi,t)=(r,\theta+\frac{2\pi}{a}, \phi+\frac{2\pi b}{a},t).
$$
Note that with this description, $m:=\gcd (a,b)>1$ is the order of the isotropy subgroup of 
$\gamma$ in the $3$-orbifold $\s^3(t_0)/\Gamma_p$. We set $a^\prime=a/m$, $b^\prime=b/m$,
and when $b^\prime\neq 0\pmod{a^\prime}$, we let $0<c<a^\prime$ be the integer satisfying 
$b^\prime c=1\pmod{a^\prime}$. In fact the case $b^\prime= 0\pmod{a^\prime}$ corresponds
to the case $a^\prime=1$, in which case we take $c=0$. With this understood, we perform the
following change of variables in the $\theta,\phi$ coordinates
$$
\theta=\theta^\prime+\frac{c}{a}\phi^\prime, \phi=\frac{1}{a^\prime}\phi^\prime.
\mbox{ where } \theta^\prime,\phi^\prime\in\R/2\pi\Z.
$$
Then the coordinates $(r,\theta^\prime,\phi^\prime,t)$ with a $\Z_m$-action given by
$$(r,\theta^\prime,\phi^\prime,t)\mapsto (r,\theta^\prime+\frac{2\pi}{m},\phi^\prime,t)$$
gives a global orbifold chart for $W/\Gamma$, a neighborhood of $\gamma$ in 
$B^4(\epsilon_0)/\Gamma_p$, on which the symplectic form $\omega_0$ takes the following form
$$
\omega_0=\frac{1}{a^\prime|\Gamma^\prime|}tdt\wedge d\phi^\prime+\alpha_0\wedge dh_0,
$$
where $\alpha_0=d\theta^\prime+\frac{c}{a}d\phi^\prime$ and $h_0=\frac{t^2r^2}{2(1+r^2)}$.
The $\s^1$-action given by translations in $\theta^\prime$ is Hamiltonian, with a Hamiltonian
function $h_0(r,\theta^\prime,\phi^\prime,t)=\frac{t^2r^2}{2(1+r^2)}$.
Let $Y(\lambda)=h_0^{-1}(\lambda)$ be the level sets. Then $Y(\lambda)$, $\lambda>0$, 
can be regarded as a principal $\s^1$-bundle over an annulus with coordinates 
$\phi^\prime,t$, and with this viewpoint $\alpha_0$ is a connection $1$-form on $Y(\lambda)$. 
Note that $Y(\lambda)$ comes with a natural trivialization. 

With these preparations,  let $\Sigma_i^\prime$ be the part of $\Sigma_i$ contained in 
$X\setminus (\bigsqcup_{p\in \Sigma^1} B^4(t_0)/\Gamma_p)$, and let $H_i$ be the subgroup of
$G$ which leaves $\Sigma_i^\prime$ invariant. Then $\Sigma_i^\prime$ is a compact Riemann 
surface with boundary where the boundary components are singular fibers in
$\bigsqcup_{p\in \Sigma^1} \s^3(t_0)/\Gamma_p$ such as $\gamma$. We can modify the arguments 
in Lemma 2.2 to de-singularize the symplectic structure $\omega$ along each $\Sigma_i^\prime$ to obtain a symplectic structure $\omega^\prime$ on 
$|X|\setminus (\bigsqcup_{p\in \Sigma^1} |B^4(t_0)/\Gamma_p|)$. 
More concretely, let $\nu_i$ be the normal bundle of $\Sigma_i^\prime$ in $X$ and let 
$\pi_i: Y_i\rightarrow \Sigma_i^\prime$ be the associated principal $\s^1$-bundle.
Note that the order $m$ of the isotropy subgroup of $\gamma$ in the $3$-orbifold 
$\s^3(t_0)/\Gamma_p$ is also the order of the isotropy group at $\Sigma_i^\prime$. With
this understood, there is a finite group $\tilde{H}_i$ with a short exact sequence 
$1\rightarrow \Z_{m}\rightarrow \tilde{H}_i\rightarrow H_i\rightarrow 1$ such that $\nu_i$
and $Y_i$ are natural $\tilde{H}_i$-bundles. We choose an $\tilde{H}_i$-invariant 
connection $1$-form $\alpha_i$ on $Y_i$, with the following condition near the boundary
of $\Sigma_i^\prime$:  in the neighborhood $W/\Gamma$ of $\gamma$, we identify the 
complement of $\Sigma_i^\prime$ with $Y_i\times (0,r_i)$ for some $r_i>0$ such that 
$Y(\lambda)=Y_i\times \{\lambda\}$, and under this identification, we require 
$\alpha_i=\alpha_0=d\theta^\prime+\frac{c}{a}d\phi^\prime$.
We let $\kappa_i$ be the $2$-form on $\Sigma_i^\prime$ such that 
$\pi_i^\ast \kappa_i=-d\alpha_i$ (note that $\kappa_i=0$ near the boundary), and consider
the following $\tilde{H}_i$-invariant symplectic form on $Y_i\times (-r_i,r_i)$ for $r_i>0$
sufficiently small:
$$
\omega_i:=\pi_i^\ast(\eta_i+\lambda \kappa_i)+\alpha_i\wedge d\lambda,
$$
where $\eta_i$ is the restriction of the symplectic form $\omega$ on $\Sigma_i^\prime$.

As in the proof of Lemma 2.2, the symplectic cutting procedure yields a symplectic structure on
a disc-bundle $D(\nu_i)$ associated to $\nu_i$ which equals $\omega_i$ under the natural 
identification of $Y_i\times (0,r_i)$ with the complement of the zero section in $D(\nu_i)$. 
We continue to denote it by $\omega_i$. Then the requirement $\alpha_i=\alpha_0=d\theta^\prime+\frac{c}{a}d\phi^\prime$
in the neighborhood $W/\Gamma$ of $\gamma$ implies easily that $\omega_i=\omega_0$ in $W/\Gamma$, where
$\omega_0$ is the original symplectic structure $\omega$, given by the following expression in the $(r,\theta^\prime,\phi^\prime,t)$
coordinates: 
$$
\omega_0=\frac{1}{a^\prime|\Gamma^\prime|}tdt\wedge d\phi^\prime+\alpha_0\wedge dh_0, \mbox{ where }
h_0=\frac{t^2r^2}{2(1+r^2)}.
$$
With this understood, the relative version of the equivariant Weinstein neighborhood theorem implies that 
$\omega_i$ gives a model for the symplectic structure $\omega$ on $X$ near each $\Sigma_i^\prime$. 
We can then use symplectic cutting as in the proof of Lemma 2.2 to extend the symplectic structure 
$\omega$ across each $\Sigma_i^\prime$. We thus obtained a symplectic structure 
$\omega^\prime$ on $|X|\setminus (\bigsqcup_{p\in \Sigma^1} |B^4(t_0)/\Gamma_p|)$, 
which is clearly $G$-invariant and has the desired properties. This finishes the proof of Lemma 2.4.

\end{proof}

\subsection{Capping-off by symplectic fillings} 
It remains to extend $\omega^\prime$ to the entire $|X|$. The key observation is that the boundary of 
$|X|\setminus (\bigsqcup_{p\in \Sigma^1} |B^4(t_0)/\Gamma_p|)$ is a concave contact boundary with respect to $\omega^\prime$, so we shall extend $\omega^\prime$ to $|X|$ by capping off the
contact boundary with an appropriate symplectic filling.  

Note that each boundary component $|\s^3(t_0)/\Gamma_p|$, $p\in\Sigma^1$,
comes with a natural Seifert fibration, which is inherited from that on the $3$-orbifold 
$\s^3(t_0)/\Gamma_p$. Furthermore, the isotropy subgroup $G_p$ acts smoothly on 
$|\s^3(t_0)/\Gamma_p|$ preserving the Seifert fibration.

\begin{lem}
For each $p\in\Sigma^1$, the boundary component $|\s^3(t_0)/\Gamma_p|$ of the symplectic
manifold $(|X|\setminus (\bigsqcup_{p\in \Sigma^1} |B^4(t_0)/\Gamma_p|), \omega^\prime)$ 
admits a $G_p$-invariant, inward-pointing Liouville vector field, such that the induced 
contact form $\alpha_p$ on $|\s^3(t_0)/\Gamma_p|$ is a $G_p$-invariant, constant multiple of
a connection $1$-form with respect to the Seifert fibration on $|\s^3(t_0)/\Gamma_p|$.
\end{lem}

\begin{proof}
First, by following through the symplectic cutting procedure in the proof of Lemma 2.4, one can easily 
show that near the boundary of $\Sigma_i^\prime$, $\omega^\prime$ is given by the expression
$$
\omega_i^\prime=\frac{1}{a^\prime|\Gamma^\prime|}tdt\wedge d\tilde{\phi}
+(d\tilde{\theta}+\frac{c}{a^\prime}d\tilde{\phi})\wedge d(\frac{t^2r^2}{2m(1+r^2)})
$$
in coordinates $(r,\tilde{\theta},\tilde{\phi},t)$,  where $\tilde{\theta},\tilde{\phi}\in\R/2\pi\Z$, 
and $\tilde{\theta},\tilde{\phi}$ are related to $\theta^\prime,\phi^\prime$ by the equations 
$\tilde{\theta}=m\theta^\prime$ and $\tilde{\phi}=\phi^\prime$.

With this understood, we consider the vector field $V:=\frac{1}{2} t\partial_t$ on 
$|B^4(\epsilon_0)/\Gamma_p|$ defined near the hypersurface $|\s^3(t_0)/\Gamma_p|$.
It is clear that $V$ is a Liouville vector field with respect to $\omega$ ($=\omega_0$) in the
complement of the singular components of the $3$-orbifold $\s^3(t_0)/\Gamma_p$. 
Near the singular components, one can easily verify that  $L_V\omega^\prime=\omega^\prime$
using the expression of $\omega_i^\prime$ given above. Finally, $V$ is clearly $G_p$-invariant,
and is inward-pointing with respect to the orientation on 
$(|X|\setminus (\bigsqcup_{p\in \Sigma^1} |B^4(t_0)/\Gamma_p|), \omega^\prime)$.

Let $\alpha_p:=i_V\omega^\prime$ be the corresponding contact form on $|\s^3(t_0)/\Gamma_p|$.
Then near the boundary of each $\Sigma_i^\prime$, the expression for $\omega^\prime_i$ implies 
that $\alpha_p$ is given by the expression
$$
\alpha_i^\prime=\frac{1}{2}t_0^2((\frac{1}{a^\prime |\Gamma^\prime|}
-\frac{cr^2}{a^\prime m(1+r^2)}) d\tilde{\phi}-\frac{r^2}{m(1+r^2)} d\tilde{\theta}).
$$
To see that $\alpha_p$ is a constant multiple of a connection $1$-form with respect to the
Seifert fibration, we note that the Seifert fibration is induced from the Hopf fibration on $\s^3(t_0)$, 
which is given by translations $\theta\mapsto \theta, \phi\mapsto \phi+\lambda$.
In terms of coordinates $\tilde{\theta},\tilde{\phi}$, it becomes 
$\tilde{\theta}\mapsto \tilde{\theta}-c\lambda, \tilde{\phi}\mapsto \tilde{\phi}+a^\prime \lambda$. 
It follows that the $\s^1$-action is generated by the vector field 
$a^\prime\partial_{\tilde{\phi}}-c\partial_{\tilde{\theta}}$, and 
$$
\alpha_i^\prime (a^\prime\partial_{\tilde{\phi}}-c\partial_{\tilde{\theta}})
=\frac{t_0^2}{2|\Gamma^\prime|}.
$$
It follows easily that $\frac{2|\Gamma^\prime|}{t^2_0}\alpha_p$ is a connection $1$-form on 
$|\s^3(t_0)/\Gamma_p|$ with respect to the Seifert fibration. Finally, 
$\alpha_p$ is clearly $G_p$-invariant. This finishes off the proof.

\end{proof}

The following lemma should be well-known to experts. For completeness and the reader's
convenience, we include a proof here.

\begin{lem}
Let $\pi: Y\rightarrow B$ be an oriented Seifert $3$-manifold, and let $\alpha_0,\alpha_1$ be connection $1$-forms on $Y$ which are also positive contact forms. Then there exists a fiber-preserving self-diffeomorphism $\psi$ such that $\psi^\ast \alpha_1=\alpha_0$. Furthermore, 
if a finite group $G$ acts on $Y$ smoothly, preserving the Seifert fibration, and both of 
$\alpha_0,\alpha_1$ are $G$-invariant, then $\psi$ can be chosen $G$-equivariant. 
\end{lem}

\begin{proof}
Let $\kappa_i$ be the $2$-form on $B$ such that $\pi^\ast \kappa_i=d\alpha_i$, $i=0,1$. 
Then since $\alpha_i$ are positive contact forms, $\kappa_i$ are symplectic forms on $B$, 
defining the same orientation. Furthermore, since $\alpha_i$ are connection $1$-forms, 
$\alpha_1-\alpha_0=\pi^\ast \eta$ for some $1$-form $\eta$ on $B$. For $B$ being $2$-dimensional, $\kappa_t:=t\kappa_1+(1-t)\kappa_0$ is symplectic for each $t\in [0,1]$. It follows that 
$\alpha_t:=t\alpha_1+(1-t)\alpha_0$ is a smooth family of connection $1$-forms 
which are also positive contact forms, and we note that 
$\frac{d}{dt}\alpha_t=\alpha_1-\alpha_0=\pi^\ast \eta$.

Let $X_t$ be the time-dependent vector field on $Y$ determined by the following conditions 
$$
X_t\in \ker \alpha_t,\;\; i_{X_t} d\alpha_t=-\pi^\ast \eta \mbox{ on $\ker \alpha_t$}.
$$
Let $Z$ be the vector field which generates the $\s^1$-action on $Y$. Then 
$$
i_{X_t} d\alpha_t(Z)=d\alpha_t (X_t,Z)=\pi^\ast \kappa_t (X_t,Z)=\kappa_t(\pi_\ast (X_t),\pi_\ast(Z))=0.
$$
Consequently, $i_{X_t} d\alpha_t=-\pi^\ast \eta$ as $1$-forms on $Y$. It follows that
$\frac{d}{dt}\alpha_t+L_{X_t}\alpha_t=0$, and the $1$-family of 
self-diffeomorphisms $\psi_t$ generated by $X_t$ obeys $\frac{d}{dt}\psi_t^\ast\alpha_t=0$, 
and $\psi_1^\ast\alpha_1=\alpha_0$ in particular. To see that $\psi_1$ preserves the $\s^1$-action, 
we note that $L_{X_t} Z=0$. This is because 
$\alpha_t(L_{X_t} Z)=-L_{X_t}\alpha_t(Z)=\pi^\ast \eta (Z)=0$ so
that $L_{X_t}Z\in\ker \alpha_t$, and on the other hand, for any $W\in \ker\alpha_t$, noting that
$\pi^\ast \kappa_t=d\alpha_t$, one has 
$$
d\alpha_t (L_{X_t} Z,W)=-L_{X_t}(d\alpha_t)(Z,W)- d\alpha_t (Z,L_{X_t}W)=\pi^\ast d\eta(Z,W)-0=0.
$$
Hence $L_{X_t} Z=0$. Finally, in the present of a fiber-preserving $G$-action, 
with $\alpha_0,\alpha_1$ being $G$-invariant, everything can be done equivariantly.
Hence the lemma. 

\end{proof}

In light of the ``uniqueness" result in Lemma 2.6, in order to complete the proof of Theorem 1.1,
we will construct an appropriate model for the contact manifolds $(|\s^3(t_0)/\Gamma_p|,\alpha_p)$
which is symplectically filled by the corresponding quotient of a $4$-ball with the standard
symplectic structure $\omega_0$. To this end, let $\Gamma^\prime_p$ be the fundamental group of 
$|\s^3(t_0)/\Gamma_p|$. Then there is a finite subgroup $\hat{\Gamma}_p$ of $U(2)$ and a short exact sequence 
$1\rightarrow \Gamma^\prime_p\rightarrow \hat{\Gamma}_p\rightarrow G_p\rightarrow 1$ such that 
$|\s^3(t_0)/\Gamma_p|$ is identified via a diffeomorphism to $\s^3/\Gamma^\prime_p$ 
and the action of $G_p$ 
on $|\s^3(t_0)/\Gamma_p|$ is induced by the action of $\hat{\Gamma}_p$ on $\s^3$. 
Furthermore, under 
this identification, the $\s^1$-action on $|\s^3(t_0)/\Gamma_p|$ defining the Seifert fibration lifts to 
a $\s^1$-action on $\s^3$, which in general takes the form 
$$\lambda \cdot (z_1,z_2)=(\lambda^mz_1,\lambda^nz_2)$$ 
for some $m,n$ with $\gcd (m,n)=1$. The key observation is that all the group actions on $\s^3$
that are involved here are complex linear, hence can be extended to $\C^2$. Furthermore,
they preserve the standard symplectic structure $\omega_0$ on $\C^2$. Note that the $\s^1$-action 
is Hamiltonian with a Hamiltonian function $h_p(z_1,z_2)=\frac{1}{2}(m r_1^2+n r_2^2)$, where 
$r_i=|z_i|$ for $i=1,2$ (note that $m,n$ depend on the point $p$). Moreover, the vector field
$L:=\frac{1}{2}(r_1\partial_{r_1}+r_2\partial_{r_2})$ is a Liouville vector field transverse to the level sets of $h_p$. We consider the ellipsoids 
$$
E(p,\delta):=\{(z_1,z_2)\in\C^2| h_p(z_1,z_2)\leq \delta\}
$$
with the standard symplectic structure $\omega_0$, and let $\alpha_0:=i_L\omega_0$ be the 
contact form on the boundary $\partial E(p,\delta)$. 

Note that $\partial E(p,\delta)$ is equivariantly diffeomorphic to the unit sphere $\s^3$ with 
respect to all the group actions involved. Consequently, $\partial E(p,\delta)/\Gamma^\prime_p$ is 
$G_p$-equivariantly diffeomorphic to $|\s^3(t_0)/\Gamma_p|$ as Seifert manifolds.
Note that $\alpha_0$ is $\hat{\Gamma}_p$-invariant and is also invariant under the $\s^1$-action. 
We let $\alpha_p^\prime$ be the descendant of $\alpha_0$ to 
$\partial E(p,\delta)/\Gamma^\prime_p$, which
is clearly $G_p$-invariant and is a constant multiple of a connection $1$-form with respect to
the Seifert fibration on $\partial E(p,\delta)/\Gamma^\prime_p$. 
By Lemma 2.6, for some choice of $\delta=\delta_p$,
there exists a $G_p$-equivariant, fiber-preserving diffeomorphic 
$\psi_p: \partial E(p,\delta_p)/\Gamma^\prime_p
\rightarrow |\s^3(t_0)/\Gamma_p|$ such that $\alpha_p^\prime=\psi^\ast_p \alpha_p$.

Note that $E(p,\delta_p)$ is $\Gamma^\prime_p$-equivariantly diffeomorphic to the $4$-ball. 
Hence the following $4$-orbifold
$$
|X|\setminus (\bigsqcup_{p\in \Sigma^1} |B^4(t_0)/\Gamma_p|)
\bigsqcup_{p\in \Sigma^1} E(p,\delta_p)/\Gamma^\prime_p,
$$
where the gluing along the boundaries is given by $\bigsqcup_{p\in\Sigma^1}\psi_p$, is diffeomorphic
to the smooth orbifold $|X|$. The former has a natural symplectic structure since each $\psi_p$ is a
contactomorphism. It is by this identification we obtained the desired symplectic $4$-orbifold
$(|X|,\omega^\prime)$. 

To finish up the proof of Theorem 1.1, we observe that since each $\psi_p$ is fiber-preserving
and $\Sigma_i\cap |\s^3(t_0)/\Gamma_p|$ is a union of fibers of the Seifert fibration on
$|\s^3(t_0)/\Gamma_p|$, the part of $\Sigma_i$ in $U_p$, which is a cone over 
$\Sigma_i\cap |\s^3(t_0)/\Gamma_p|$, can be identified with a cone in 
$E(p,\delta_p)/\Gamma^\prime_p$ over 
$\psi_p^{-1}(\Sigma_i\cap |\s^3(t_0)/\Gamma_p|)\subset \partial E(p,\delta_p)/\Gamma^\prime_p$.
We should point out that with respect to the smooth structure of $E(p,\delta_p)/\Gamma^\prime_p$, 
it can be singular at the origin $0\in E(p,\delta_p)$ because the fiber of the Seifert fibration on
$\partial E(p,\delta_p)$ is in general a $(m,n)$-torus knot. (Examples in \cite{Dun} show that 
this can indeed occur.) However, with respect to the standard complex structure on
$E(p,\delta_p)$, the cone is given by a holomorphic curve. This completes the proof of Theorem 1.1.

\vspace{2mm}

We end this section with a proof of Corollary 1.3.

\vspace{2mm}

\noindent{\bf Proof of Corollary 1.3}

\vspace{2mm}

We apply Theorem 1.1 to the symplectic $4$-orbifold $\hat{M}/\Gamma$, which is equipped 
with a symplectic $G/\Gamma$-action. The underlying space $|\hat{M}/\Gamma|$ is naturally identified with
$M$ and the singular set $\Sigma$ of $\hat{M}/\Gamma$ is identified with the branch locus $B$. 
With this understood, it follows that $M$ admits a $G/\Gamma$-invariant symplectic structure with 
respect to which $B$ is a symplectic surface. Note that in the present case, $\Sigma^1$ is empty, 
and so is $U$.

For the converse suppose $M$ admits a $G/\Gamma$-invariant symplectic structure with respect to 
which $B$ is a symplectic surface. Then if we let $\hat{B}$ be the pre-image of $B$ in $\hat{M}$,
then the $G/\Gamma$-invariant symplectic structure on $M$ lifts to a $G$-invariant symplectic
structure on $\hat{M}\setminus \hat{B}$. To see that it extends across $\hat{B}$ to a symplectic
structure on $\hat{M}$, we simply observe that in the proof of Lemma 2.2, if we are given a symplectic 
structure on the disc bundle $D(\nu^\prime)$, we can lift it to a symplectic structure on $D^\ast(\nu,r)$, 
and then the same symplectic cutting procedure will allow us to extend the symplectic structure on 
$D^\ast(\nu,r)$ across the zero section to a symplectic structure on the entire disc bundle 
$D(\nu,r)$. (Compare also \cite{MRT}, Proposition 7, for a different construction to this effect.)
This completes the proof of Corollary 1.3.

\section{Symplectic resolution and its various properties}
This section is devoted to the construction of symplectic resolution and its various properties. 
In particular, we present a proof of Theorem 1.5. The section also contains the proofs of two propositions, one concerning the canonical class of the symplectic resolution, the other concerning symplectic equivariant blowing-down.

\subsection{Construction of symplectic resolution}
As we mentioned in the introduction, the resolution $\tilde{X}$ in Theorem 1.5 is simply the symplectic resolution of the symplectic $4$-orbifold $(|X|,\omega^\prime)$. For the specific purpose of applications in finite group actions,  we shall adopt the method in \cite{CFM} for its construction.

\vspace{2mm}

\noindent{\bf Proof of Theorem 1.5}

\vspace{2mm}

First, we construct the resolution $\tilde{X}$, together with the map $\pi: \tilde{X}\rightarrow X$. 
Since this is a local operation, it suffices to focus at a point $p\in \Sigma^0\sqcup\Sigma^1$,
where when $p\in \Sigma^1$, we assume it is a singular point of $(|X|,\omega^\prime)$.
For simplicity, we denote by $\tilde{\Sigma}^1$ the subset of $\Sigma^1$ which consists of
singular points of $(|X|,\omega^\prime)$. 

By the equivariant version of Darboux theorem, there is a neighborhood $U_p$ of $p$ in the orbifold $(|X|,\omega^\prime)$, which is modeled by $(B^4,\omega_0)/\Gamma^\prime_p$,
where $B^4\subset \C^2$ is a $4$-ball, $\omega_0$ is the standard symplectic structure, and
$\Gamma_p^\prime$ is the isotropy group at $p$, acting on $\C^2$ as a subgroup of $U(2)$.
If $p\in\tilde{\Sigma}^1$, we shall assume $B^4$ is contained in the ellipsoid $E(p,\delta_p)$ from
the proof of Theorem 1.1. Note that the $\omega^\prime$-compatible complex structure in
Theorem 1.1, with respect to which $\Sigma_i\cap U$ is given by a holomorphic curve, is
simply the standard complex structure on $B^4$ in the case of $p\in\tilde{\Sigma}^1$. 
We identify $U_p=B^4/\Gamma_p^\prime$
as a subset of the affine algebraic variety $Z:=\C^2/\Gamma_p^\prime$.

Let $\pi_Z: \tilde{Z}\rightarrow Z$ be the minimal algebraic resolution of the isolated singularity
of $Z$. Setting $\tilde{U}_p:=\pi_Z^{-1}(U_p)$, where $U_p=B^4/\Gamma_p^\prime
\subset Z$, we glue each $\tilde{U}_p$ to $|X|\setminus (\Sigma^0\sqcup\Sigma^1)$ by
identifying $\tilde{U}_p\setminus \pi_Z^{-1}(p)$ with $U_p\setminus\{p\}$ via $\pi_Z$, which
gives us the desired resolution $\tilde{X}$ together with the map $\pi: \tilde{X}\rightarrow X$. 
The claims in Theorem 1.5(1) follow easily. We remark that since we require that
the algebraic resolution $\pi_Z: \tilde{Z}\rightarrow Z$ be minimal, which is unique (cf. \cite{Lau}), 
the diffeomorphism type of $\tilde{X}$ is canonically determined by the smooth orbifold $X$. 

Secondly, fixing any such a choice of $U_p$, $p\in \Sigma^0\sqcup\tilde{\Sigma}^1$, 
we construct a symplectic structure $\tilde{\omega}$ on $\tilde{X}$ as follows. First, 
we note that $\tilde{Z}$ is a quasi-projective variety, and hence there is a K\"{a}hler form 
$\Omega$ on  $\tilde{Z}$. Consider for each $p$ the subset 
 $A:=\pi_Z^{-1}((\frac{2}{3} \overline{B}^4\setminus \frac{1}{3} B^4)/\Gamma^\prime_p)$
 in $\tilde{U}_p$. Since it is homotopic to $\s^3/\Gamma_p^\prime$, a rational homology sphere,
 there is a $1$-form $\gamma\in \Omega^1(A)$ such that 
 $\Omega-\pi_Z^\ast \omega^\prime=d\gamma$ on $A$.
 We let $\rho$ be a cut-off function which equals zero in 
 $\tilde{U}_p\setminus \pi_Z^{-1}(\frac{2}{3} B^4/\Gamma_p^\prime)$ and equals one
 in $\pi_Z^{-1}(\frac{1}{3} B^4/\Gamma_p^\prime)$. Then we define a closed $2$-form 
 $\tilde{\omega}$ on $\tilde{U}_p$, where
 $$
\tilde{ \omega}:=\pi_Z^\ast \omega^\prime+ \epsilon d(\rho\gamma)
 $$
 for some $\epsilon>0$ sufficiently small. Note that on 
 $\tilde{U}_p\setminus  \pi_Z^{-1}(\frac{2}{3} B^4/\Gamma_p^\prime)$, $\tilde{\omega}=
 \pi_Z^\ast \omega^\prime$ which is symplectic, 
 and on $\pi_Z^{-1}(\frac{1}{3} B^4/\Gamma_p^\prime)$,
 $\tilde{\omega}=(1-\epsilon)\pi_Z^\ast\omega^\prime+\epsilon\Omega$, which is a
 K\"{a}hler form, and finally on $A$, $\tilde{\omega}$ is symplectic as long as we choose
 $\epsilon>0$ sufficiently small since $A$ is compact. It is clear that $\tilde{\omega}$ extends
 to the whole $\tilde{X}$, giving the desired symplectic structure. Let $U$ be the $G$-invariant neighborhood of $\Sigma^0\sqcup\Sigma^1$, where
 $$
 U:=\bigsqcup_{p\in \Sigma^0} U_p\bigsqcup_{p\in\Sigma^1} E(p,\delta_p)/\Gamma_p,
 $$
it follows easily that  
$\pi: (\tilde{X}\setminus \pi^{-1}(U\cup\Sigma^\ast),\tilde{\omega})\rightarrow (X\setminus (U\cup \Sigma^\ast),\omega)$ 
is a symplectomorphism. Note that $U$ can be taken arbitrarily small. To see that 
$\pi^{-1}(\Sigma^\ast)$ is a symplectic surface, we only need to observe that on $A$ the complex structure is $\tilde{\omega}$-tame and $A\cap \pi^{-1}(\Sigma^\ast)$ is holomorphic. 
Furthermore, it is clear that $\tilde{\omega}=\pi^\ast \omega$ on $\pi^{-1}(\Sigma^\ast\setminus U)$ as area forms. Finally, we take 
$U^\prime:=\bigsqcup_{p\in\tilde{\Sigma}^1} \pi_Z^{-1}(\frac{1}{3} B^4/\Gamma_p^\prime)$. 
Then $U^\prime\cap \pi^{-1}(\Sigma)$ is given by holomorphic curves. The rest of the claims in Theorem 1.5(2) are obvious.
  
 Finally, we prove the claims in Theorem 1.5(3). First, it is well-known that the algebraic
 resolution $\pi_Z:\tilde{Z}\rightarrow Z$ can be carried out equivariantly (cf. \cite{V}). It follows easily
 that if the orbifold $(X,\omega)$ admits a symplectic $G$-action by a finite group $G$, then
 the constructions of the resolution $\pi: \tilde{X}\rightarrow X$ and the symplectic structure
 $\tilde{\omega}$ can be done $G$-equivariantly. 

It remains to compare the resolutions of the symplectic orbifolds $X/G$ and $\tilde{X}/G$.
To this end, using the same notations from the construction of $\tilde{X}$, we consider the 
$G$-invariant decompositions of the underlying spaces of $X$ and $\tilde{X}$:
 $|X|$ as the union of $|X|\setminus (\bigsqcup_{p\in \Sigma^0\sqcup \tilde{\Sigma}^1} U_p)$
 and $\bigsqcup_{p\in \Sigma^0\sqcup \tilde{\Sigma}^1} U_p$, and 
 $\tilde{X}$ as the union of $|X|\setminus (\bigsqcup_{p\in \Sigma^0\sqcup \tilde{\Sigma}^1} U_p)$
 and $\bigsqcup_{p\in \Sigma^0\sqcup \tilde{\Sigma}^1} \tilde{U}_p$. Noting that the underlying space 
 of $X/G$ is simply the underlying space of the quotient of $|X|$ with the induced $G$-action,
 this gives rise to the corresponding decompositions of the underlying spaces of $X/G$ and 
 $\tilde{X}/G$:
 $$
 |X/G|=(|X|\setminus (\bigsqcup_{p\in \Sigma^0\sqcup \tilde{\Sigma}^1} U_p))/G
 \bigsqcup_{p\in (\Sigma^0\sqcup \tilde{\Sigma}^1)/G} |U_p/G_p|
 $$
 and
 $$
 |\tilde{X}/G|=(|X|\setminus (\bigsqcup_{p\in \Sigma^0\sqcup \tilde{\Sigma}^1} U_p))/G
 \bigsqcup_{p\in (\Sigma^0\sqcup \tilde{\Sigma}^1)/G} |\tilde{U}_p/G_p|,
$$
where $G_p$ is the subgroup of $G$ fixing $p$, and $p\in (\Sigma^0\sqcup \tilde{\Sigma}^1)/G$ means that $p$ is running over a set of representatives of the quotient set 
$(\Sigma^0\sqcup \tilde{\Sigma}^1)/G$ in $\Sigma^0\sqcup \tilde{\Sigma}^1$.

With the preceding understood, let $\pi_V:V\rightarrow X/G$ and $\pi_W: W\rightarrow \tilde{X}/G$
denote the corresponding resolutions. It follows easily that the difference between $W$ and $V$ 
occurs at $\pi_V^{-1}(U_p/G_p)$ and $\pi_W^{-1}(\tilde{U}_p/G_p)$, where 
$p\in (\Sigma^0\sqcup \tilde{\Sigma}^1)/G$.
With this understood, the following claim finishes off the proof of Theorem 1.5(3).

\vspace{2mm}

{\bf Claim:}{\it \hspace{1mm} $\pi_W^{-1}(\tilde{U}_p/G_p)$ is either diffeomorphic to
$\pi_V^{-1}(U_p/G_p)$, or can be reduced to a manifold diffeomorphic to it by successively 
blowing down symplectic $(-1)$-spheres. 
}

\vspace{2mm}

{\bf Proof of Claim:}
We begin by noting that $|U_p/G_p|$ is a complex analytic space with a unique, isolated 
quotient singularity at $p$. Let $\widetilde{|U_p/G_p|}$ be its minimal resolution. On the other
hand, $|\tilde{U}_p/G_p|$ is a complex analytic space with isolated quotient
singularities, all contained in $\pi_Z^{-1}(p)/G_p$. Let $\widetilde{|\tilde{U}_p/G_p|}$ denote
the minimal resolution of $|\tilde{U}_p/G_p|$. It is easily seen that 
$\widetilde{|\tilde{U}_p/G_p|}$ is some resolution of the unique singularity of $|U_p/G_p|$.

With the preceding understood, we note, from the construction of the resolutions $V$ and $W$, 
that $\pi_V^{-1}(U_p/G_p)=\widetilde{|U_p/G_p|}$ and $\pi_W^{-1}(\tilde{U}_p/G_p)$ is
diffeomorphic to $\widetilde{|\tilde{U}_p/G_p|}$. Moreover, regarding 
$\widetilde{|\tilde{U}_p/G_p|}$ as a resolution of the unique singularity in $|U_p/G_p|$,
the exceptional set in $\widetilde{|\tilde{U}_p/G_p|}$ corresponds to a configuration of embedded 
symplectic two-spheres in $\pi_W^{-1}(\tilde{U}_p/G_p)$ intersecting transversely. Furthermore, 
there is an almost complex structure $J$ compatible with the symplectic structure such that each 
symplectic two-sphere in the configuration is $J$-holomorphic. 

By the work of Artin \cite{Artin}, the exceptional set in $\widetilde{|\tilde{U}_p/G_p|}$
has the following properties: any two distinct components are either disjoint or intersect 
at a single point, and no three distinct components intersect in one point. Furthermore, 
the dual graph is a tree. Clearly, the corresponding configuration of symplectic two-spheres
in $\pi_W^{-1}(\tilde{U}_p/G_p)$ also has these properties. 

Now we recall the following fact: $\widetilde{|\tilde{U}_p/G_p|}$ can be reduced to the minimal
resolution $\widetilde{|U_p/G_p|}$ by successively blowing down holomorphic $(-1)$-spheres 
(cf. \cite{Lau}). The following lemma, Lemma 3.1,  shows that the corresponding successive 
blowing-downs for 
$\pi_W^{-1}(\tilde{U}_p/G_p)$ can be done symplectically, i.e., with the holomorphic $(-1)$-spheres 
replaced by symplectic $(-1)$-spheres in each step. With this understood, it follows easily that 
either $\pi_W^{-1}(\tilde{U}_p/G_p)$ is diffeomorphic to $\pi_V^{-1}(U_p/G_p)$ or can be reduced 
to a manifold diffeomorphic to $\pi_V^{-1}(U_p/G_p)$. Hence the claim. 

\begin{lem}
Let $S$ be a symplectic $(-1)$-sphere in a symplectic $4$-manifold $(M,\omega)$, and let $\{C_i\}$ be
a finite collection of symplectic surfaces in $(M,\omega)$. Suppose there is an $\omega$-compatible almost complex structure 
$J$ such that $S$ and $C_i$ are $J$-holomorphic, and furthermore, near $S$ the symplectic surfaces
$C_i$ are embedded, disjoint, and each $C_i$ intersects $S$ transversely. Then 
\begin{itemize}
\item [{(1)}] there is a neighborhood $U$ of $S$,  such that for each $i$, there is an embedded 
symplectic surface $\tilde{C}_i$, which is isotopic to $C_i$ through a symplectic isotopy 
supported in $U$, 
\item [{(2)}] there is a smaller neighborhood $V\subset U$ such that, if we let $(M^\prime,\omega^\prime)$ be the 
symplectic $4$-manifold obtained by removing $V$ and gluing back a standard symplectic $4$-ball, then
each symplectic surface $\tilde{C}_i\setminus V$ can be naturally extended across the $4$-ball to a closed symplectic surface $C_i^\prime$, such that the surfaces $C_i^\prime$ intersect transversely 
at the origin of the $4$-ball, and 
\item [{(3)}] there is an $\omega^\prime$-compatible almost complex structure $J^\prime$,
agreeing with $J$ outside $U$, such that $C_i^\prime$ is $J^\prime$-holomorphic. 
\end{itemize}
\end{lem}

\begin{proof}
First of all, one can isotop each $C_i$ near the intersection points with $S$, which can be made in an arbitrarily small
neighborhood of the intersection points, such that the new symplectic surface is $\omega$-orthogonal to $S$ (cf. \cite{Gom}).
With this understood, let $U$ be a neighborhood of $S$ in which the symplectic form $\omega$ is modeled by a standard 
symplectic structure on the disc bundle of $S$. This is possible by the Weinstein neighborhood theorem. Now observe that
in the standard model the fiber disc is $\omega$-orthogonal to $S$. By further deforming each $C_i$ inside $U$, we can 
arrange so that near the intersection point it coincides with the fiber disc. This final new symplectic surface is our $\tilde{C}_i$.
Now we take a sufficiently small neighborhood $V\subset U$ of $S$ such that inside $V$, each $\tilde{C}_i$ is given by
the fiber disc. Then clearly, after removing $V$ and gluing back a standard symplectic $4$-ball, each $\tilde{C}_i\setminus V$ 
can be extended across the $4$-ball by gluing a standard complex disc to the boundary circle of $\tilde{C}_i\setminus V$.
The resulting closed symplectic surfaces $C_i^\prime$ intersect transversely at the origin of the $4$-ball. We define the
$\omega^\prime$-compatible almost complex structure $J^\prime$ to be $J$ outside the neighborhood $U$, to be the 
standard complex structure inside the $4$-ball, and to be some $\omega^\prime$-compatible almost complex structure in
$U\setminus V$ such that each $\tilde{C}_i$ is pseudo-holomorphic. The last assertion is possible because inside 
$U\setminus V$ the symplectic surfaces $\tilde{C}_i$ are embedded and disjoint. It is clear that each $C_i^\prime$ is 
$J^\prime$-holomorphic. This finishes off the lemma.
 
\end{proof}

The proof of Theorem 1.5 is complete.

\subsection{The canonical class of the symplectic resolution}
As we pointed out in Remark 1.2, the orbifold canonical line bundle $K_{\omega^\prime}$ in
Theorem 1.1 is uniquely determined up to isomorphisms.
To see this, we first note that on $|X|\setminus U$, $\omega^\prime$ is uniquely determined by 
the original symplectic structure $\omega$ because $\omega=\omega^\prime$ on
$|X|\setminus (\Sigma^\ast \cup U)$ and $|X|\setminus (\Sigma^\ast \cup U)$ is dense in
$|X|\setminus U$. Consequently, the restriction of $K_{\omega^\prime}$ over $|X|\setminus U$
is uniquely determined by $(X,\omega)$. On the other hand, we observe that up to isomorphism, 
there is only one way to extend $K_{\omega^\prime}|_{|X|\setminus U}$ over $U$. Hence 
$K_{\omega^\prime}$ is uniquely determined up to isomorphism. We also pointed out in 
Remark 1.6 that the canonical line bundle $K_{\tilde{\omega}}$ is uniquely determined up to isomorphism. This is because $K_{\tilde{\omega}}$ is completely determined by 
$K_{\omega^\prime}$ and the singularities of the smooth orbifold $|X|$ (recall that we have 
used the unique, minimal resolutions of the singularities of $|X|$ in the construction of $\tilde{X}$).

In the following proposition, we give an expression of $c_1(K_{\tilde{X}})$ (as a class in 
$H^2(\tilde{X},\Q)$) in terms of the orbifold canonical class $c_1(K_X)$ and the singularities 
of $(|X|,\omega^\prime)$.

\begin{prop}
Let $\{p_j\}$ be the set of singular points of $(|X|,\omega^\prime)$, and for each $p_j$, let 
$\{E_{j,k}|k\in I_j\}$ be the exceptional set in the minimal resolution of $p_j$. Then there are
$a_{j,k}\in \Q$, $a_{j,k}\leq 0$, such that
$$
c_1(K_{\tilde{X}})=\pi^\ast c_1(K_{|X|})+ \sum_{j} \sum_{k\in I_j} a_{j,k} E_{j,k}.
$$
On the other hand, let $\{\Sigma_i\}$ be the set of compactified connected components of 
$\Sigma^\ast$, and for each $i$, let $m_i$ be the order of the isotropy groups along $\Sigma_i$.
Then $c_1(K_{|X|})$ and $c_1(K_X)$ are related by the following equation
$$
c_1(K_{|X|})=c_1(K_{X}) +\sum_i \frac{1-m_i}{m_i}PD(\Sigma_i),
$$
where $PD(\Sigma_i)$ stands for the Poincar\'{e} dual of the symplectic surface $\Sigma_i$
in $|X|$.
\end{prop}

\begin{proof}
For each $j$, we fix a regular neighborhood $V_j$ of $p_j$ and let $\tilde{V}_j$ be the minimal
resolution of $V_j$ at $p_j$. Then by the Mayer-Vietoris theorem (with $\Q$-coefficients), 
$$
c_1(K_{\tilde{X}})=c_1(K_{|X|\setminus \sqcup_j V_j})+\sum_j c_1(K_{\tilde{V}_j}).
$$
On the other hand, note that $c_1(K_{|X|\setminus \sqcup_j V_j})=\pi^\ast c_1(K_{|X|})$,
$c_1(K_{\tilde{V}_j})=\sum_{k\in I_j} a_{j,k} E_{j,k}$ for some $a_{j,k}\in\Q$, where $a_{j,k}\leq 0$
as $\tilde{V}_j$ is the minimal resolution. The formula 
 $$
c_1(K_{\tilde{X}})=\pi^\ast c_1(K_{|X|})+ \sum_{j} \sum_{k\in I_j} a_{j,k} E_{j,k}.
$$
follows immediately.

The main part of the proof is concerned with the equation relating $c_1(K_{|X|})$ and $c_1(K_X)$.
To this end, we introduce the following notations: $Z:=X\setminus U$, $|Z|:=|X|\setminus U$,
and $\Sigma_i^\prime:=\Sigma_i\setminus U$, where $U$ is a neighborhood of $\Sigma^1$ introduced in Theorem 1.1. With this understood, note that the symplectic structures $\omega$ 
and $\omega^\prime$ agree in the complement of 
$\Sigma^\ast\setminus U=\bigsqcup_i \Sigma_i^\prime$. This implies that for any fixed regular neighborhood $U^\prime$ of $\Sigma^\ast\setminus U$ in $|Z|$, we can choose an $\omega$-compatible 
almost complex structure $J$ on $Z$ and an $\omega^\prime$-compatible almost complex 
structure $J^\prime$ on $|Z|$ such that $J=J^\prime$ in the complement of $U^\prime$. 
This in particular yields an identification of $K_Z$ and $K_{|Z|}$ in the complement of $U^\prime$. 
On the other hand, note that there is a continuous orbifold map $\lambda: Z\rightarrow |Z|$ 
which induces the identity map between the underlying spaces (cf. \cite{C01}). Furthermore, 
$\lambda$ is smooth in the complement of $\Sigma^\ast\setminus U$. 
Now if we denote by $K_{|Z|}^\ast$
the dual of $K_{|Z|}$, we obtain the following bundle isomorphisms 
$$
K_{Z}=\lambda^\ast K_{|Z|}\otimes \lambda^\ast K_{|Z|}^\ast\otimes K_Z
=\lambda^\ast K_{|Z|}\otimes \text{Hom}(\lambda^\ast K_{|Z|},K_Z),
$$
where $\text{Hom}(\lambda^\ast K_{|Z|},K_Z)$ is the bundle of endomorphisms from 
$\lambda^\ast K_{|Z|}$ to $K_Z$. With this understood, note that the identification of 
$K_Z$ and $K_{|Z|}$ in the complement of $U^\prime$ defines a canonical non-zero 
section $\sigma$ of $\text{Hom}(\lambda^\ast K_{|Z|},K_Z)$ in the complement 
of $U^\prime$. It follows that the first Chern class of $\text{Hom}(\lambda^\ast K_{|Z|},K_Z)$ 
should be given by a linear combination of the Poincar\'{e} duals of
$\Sigma_i^\prime$. (Note that each $\Sigma_i^\prime$ defines a class in 
$H_2(|Z|,\partial |Z|,\Q)$, so its Poincar\'{e} dual lies in $H^2(|Z|,\Q)$.)

In order to compute the first Chern class of $\text{Hom}(\lambda^\ast K_{|Z|},K_Z)$, 
we let $U^\prime_i$ be the component of $U^\prime$ containing $\Sigma_i^\prime$, 
which is taken to be a disc bundle of the normal bundle $\pi_i^\prime: \nu_i^\prime\rightarrow \Sigma_i^\prime$ in $|Z|$. Let  $U_i$ be the global orbifold chart of $Z$ over $U_i^\prime$, 
which is also a disc bundle over $\Sigma_i^\prime$ associated to the normal bundle 
$\pi_i: \nu_i\rightarrow \Sigma_i^\prime$ of $\Sigma_i^\prime$ in $Z$ (see the proof of Lemma 2.4). The $\Z_{m_i}$-action on $U_i$ is given by the complex multiplication on the fibers. 
Let $\lambda_i: U_i\rightarrow U_i^\prime=U_i/\Z_{m_i}$ be the quotient map. Then with 
this understood, the bundle $\text{Hom}(\lambda^\ast K_{|Z|},K_Z)$ over $U_i^\prime$ is 
given by the following $\Z_{m_i}$-equivariant bundle over $U_i$:
$$
\text{Hom}(\lambda_i^\ast \circ (\pi_i^\prime)^\ast (\nu_i^\prime\otimes T\Sigma_i^\prime)^\ast, 
\pi_i^\ast (\nu_i\otimes T\Sigma_i^\prime)^\ast)=
\text{Hom}(\lambda_i^\ast (\pi_i^\prime)^\ast (\nu_i^\prime)^\ast\otimes 
\pi_i^\ast T^\ast\Sigma_i^\prime, \pi_i^\ast \nu_i^\ast\otimes \pi_i^\ast T^\ast\Sigma_i^\prime)
$$ 
The canonical non-zero section $\sigma$ of $\text{Hom}(\lambda^\ast K_{|Z|},K_Z)$ in the complement of $U^\prime$ determines an $\Z_{m_i}$-equivariant non-zero section $\sigma_i$ of 
$$
\text{Hom}(\lambda_i^\ast (\pi_i^\prime)^\ast (\nu_i^\prime)^\ast\otimes 
\pi_i^\ast T^\ast\Sigma_i^\prime, \pi_i^\ast \nu_i^\ast\otimes \pi_i^\ast T^\ast\Sigma_i^\prime)
$$
on $\overline{U_i}\setminus U_i$. Now observe that on the factor $\pi_i^\ast T^\ast\Sigma_i^\prime$, 
$\sigma_i$ is given by the identity map. In order to understand $\sigma_i$ between the factors
$\lambda_i^\ast (\pi_i^\prime)^\ast (\nu_i^\prime)^\ast$ and $\pi_i^\ast \nu_i^\ast$, we let 
$D_p\subset U_i$ be the disc which is the fiber of $\pi_i$ over $p\in \Sigma_i^\prime$, and let
$D_p^\prime\subset U_i^\prime$ be the disc which is the fiber of $\pi_i^\prime$. Denote by
$z$ and $w$ the complex coordinates on $D_p$ and $D_p^\prime$ respectively. Then over $D_p$,
$\lambda_i^\ast (\pi_i^\prime)^\ast (\nu_i^\prime)^\ast$, $\pi_i^\ast \nu_i^\ast$ are trivialized by
$dw$ and $dz$. On the other hand, since $J=J^\prime$ in the complement of $U^\prime$, 
$dw$ and $dz$ are related on $\partial D_p$ by the equation $dw=z^{m_i-1} dz$. Now by
the fact that the non-zero section $\sigma$ is defined by the identification of 
$K_Z$ and $K_{|Z|}$ in the complement of $U^\prime$, it follows easily that $\sigma_i$
on $\partial D_p$ is given by the map which assigns each $z\in \partial D_p$ the automorphism 
$\sigma_i(z)$ of $\C$ given by the multiplication by $z^{m_i-1}$. In other words, 
$\sigma_i:\partial D_p\times \C\rightarrow \partial D_p\times \C$ is given by the formula 
$$
\sigma_i: (z,u)\mapsto (z,z^{m_i-1}u),
$$
which can be naturally extended to a map $D_p\times \C\rightarrow D_p\times \C$ given by
the same formula above. It follows easily that the first Chern class of 
$\text{Hom}(\lambda^\ast K_{|Z|},K_Z)$ equals the Poincar\'{e} dual of $\sum_i \frac{m_i-1}{m_i}\Sigma_i^\prime$. Now recall that $\lambda: Z\rightarrow |Z|$ induces the identity map on 
$H^2(|Z|;\Q)$, it follows immediately that 
$$
c_1(K_Z)=c_1(K_{|Z|}) +\sum_i \frac{m_i-1}{m_i}PD(\Sigma_i^\prime).
$$

Finally, we observe that $U$ is a disjoint union of $\Q$-homology balls. Hence there is a
natural identification between $H^2(|X|,\Q)$ and $H^2(|Z|,\Q)$. It is easy to see that under
this identification, $c_1(K_X)=c_1(K_Z)$, $c_1(K_{|X|})=c_1(K_{|Z|})$, and 
$PD(\Sigma_i)=PD(\Sigma_i^\prime)$ for each $i$. This implies immediately 
$$
c_1(K_{|X|})=c_1(K_{X}) +\sum_i \frac{1-m_i}{m_i}PD(\Sigma_i),
$$
and the proof of the proposition is complete.

\end{proof}

\subsection{Symplectic resolution and equivariant blowing down}
The symplectic blowing down operation (cf. \cite{McS}) can be easily extended to the 
equivariant setting. More concretely, let $\tilde{M}$ be a symplectic $4$-manifold 
equipped with a finite symplectic $G$-action. Suppose there exists a $G$-invariant set of disjoint symplectic $(-1)$-spheres in $\tilde{M}$. Then blowing down $\tilde{M}$ 
along the $(-1)$-spheres gives rise to a symplectic $4$-manifold $M$, which can be arranged 
so that the $G$-action descends to a symplectic $G$-action on $M$. The symplectic $G$-manifold 
$\tilde{M}$ is called {\it minimal} if no such set of $(-1)$-spheres exists. 
We refer the reader to \cite{C2} for more discussions on this topic. 

Since it is technically more convenient to work with minimal symplectic $G$-manifolds, one naturally asks how the resolutions of the quotient orbifolds are related after performing an equivariant symplectic blowing-down. 

Recall that for a symplectic $4$-manifold $M$ equipped with a finite symplectic 
$G$-action, the resolution of the quotient orbifold $M/G$ is denoted by $M_G$.
The arguments in the proof of Theorem 1.5(3) can be easily extended to give a proof of the 
following proposition.

\begin{prop}
Let $\tilde{M}$ be a symplectic $4$-manifold equipped with a finite symplectic $G$-action. 
Suppose $M$ is a $G$-equivariant blow-down of $\tilde{M}$. Then the resolutions $\tilde{M}_G$ 
and $M_G$ are either diffeomorphic, or $\tilde{M}_G$ can be reduced to $M_G$ by successively blowing
down symplectic $(-1)$-spheres. In particular, $\kappa^s(\tilde{M}_G)=\kappa^s(M_G)$.
\end{prop}

\begin{proof}
Note that as in the proof of Theorem 1.5(3), the issue is local in nature. For this reason it suffices to consider a symplectic $(-1)$-sphere
$C$ in $\tilde{M}$ which is blown down to a point $p$ in $M$. Let $G_C$ be the subgroup of $G$
which leaves $C$ invariant. The key observation is that the $G_C$-action in a neighborhood
of $C$ is smoothly equivalent to a holomorphic action. To see this, note that $C$ has
a $G_C$-invariant regular neighborhood whose boundary is $\s^3$. This implies that
$G_C$ acts smoothly and effectively on $\s^3$, which preserves the Hopf fibration.
Any such an action is equivalent to a linear action preserving the Hopf fibration, hence 
it is by a subgroup of $U(2)$. Now we fix such a complex linear action of $G_C$ on
$\C^2$. Equivariantly blowing up at the origin of $\C^2$, we obtain a holomorphic model 
for the $G_C$-action near $C$. 

With the preceding understood, let $U_C$ be a $G_C$-invariant neighborhood of $C$ 
in $\tilde{M}$. Note that the holomorphic model in the previous paragraph supports a
$G_C$-invariant K\"{a}hler form. Hence by the equivariant Weinstein neighborhood theorem,
we may assume that there is a complex structure on $U_C$, compatible with the symplectic
structure on $\tilde{M}$, such that $C$ is holomorphic and the $G_C$-action is holomorphic. 
Let $\pi_C: U_C\rightarrow U_p$ be the holomorphic map contracting $C$ to a point. 
Then there is an induced holomorphic $G_C$-action on $U_p$, and the equivariant 
symplectic blowing down operation from $\tilde{M}$ to $M$ is locally near $C$ smoothly 
equivalent to $\pi_C: U_C\rightarrow U_p$. This said, if we let $\tilde{\pi}: \tilde{M}_G
\rightarrow \tilde{M}/G$ and $\pi: M_G\rightarrow M/G$ be the corresponding resolutions,
then $\pi^{-1}(U_p/G_C)$ is diffeomorphic to the minimal resolution of the unique isolated 
singularity of the complex analytic space $|U_p/G_C|$, and $\tilde{\pi}^{-1}(U_C/G_C)$ is 
diffeomorphic to some resolution of the singularity of $|U_p/G_C|$. By the same arguments 
as in the proof of Theorem 1.5(3), with the help of Lemma 3.1, $\tilde{\pi}^{-1}(U_C/G_C)$ 
can be reduced to a manifold diffeomorphic to $\pi^{-1}(U_p/G_C)$ by successively blowing 
down symplectic $(-1)$-spheres. Proposition 3.3 follows easily from this local consideration. 

\end{proof}

\section{Symplectic Kodaira dimension of $M_G$}

\subsection{The case of $\kappa^s(M)=-\infty$ or $0$}

We begin with the following lemma.

\begin{lem}
Let $(M,\omega)$ be a symplectic $4$-manifold equipped with a finite symplectic $G$-action.
Suppose one of the following conditions is satisfied:
\begin{itemize}
\item $c_1(K_M)\cdot [\omega]<0$, or
\item $M$ has torsion canonical class, and the singular set of $M/G$ either contains a
$2$-dimensional component, or contains an isolated non-Du Val singularity.
\end{itemize}
Then one can choose the symplectic structure $\tilde{\omega}$ on $M_G$ such that
$c_1(K_{M_G})\cdot [\tilde{\omega}]<0$. In particular, $M_G$ is rational or ruled.
\end{lem}

\begin{proof}
Let $X=M/G$ be the quotient orbifold. For simplicity, we continue to denote by $\omega$ the 
symplectic structure on $X$. We let $\{\Sigma_i\}$ be the $2$-dimensional components of the
singular set of $X$, regarded as symplectic surfaces in the orbifold $(|X|,\omega^\prime)$,
and let $\{p_j\}$ be the set of singular points of $(|X|,\omega^\prime)$.  With this understood, 
let $m_i>1$ be the order of the isotropy groups along $\Sigma_i$,  and for each $j$, 
let $\{E_{j,k}|k\in I_j\}$ be the set of exceptional divisors in the minimal resolution of the singular 
point $p_j$.

By Proposition 3.2, we have
$$
c_1(K_{|X|})=c_1(K_X)+\sum_i \frac{1-m_i}{m_i} PD(\Sigma_i),
$$
and
$$
c_1(K_{M_G})=\pi^\ast c_1(K_{|X|}) +\sum_j \sum_{k\in I_j} a_{j,k} E_{j,k}
$$
where $a_{j,k}\in \Q$ and $a_{j,k}\leq 0$. With this understood, we compute 
$c_1(K_{M_G})\cdot [\tilde{\omega}]$. 

To this end, we fix a closed $2$-form $\eta$ representing $c_1(K_{|X|})$. Then recall that in the definition of $\tilde{\omega}$ in the proof of Theorem 1.5, 
$\tilde{ \omega}:=\pi_Z^\ast \omega^\prime+ \epsilon d(\rho\gamma)$ on $\tilde{U}_p$.
It follows easily that 
$$
|\pi^\ast c_1(K_{|X|}) \cdot [\tilde{\omega}] - c_1(K_{|X|})\cdot [\omega^\prime]|\leq 
\epsilon \cdot |\eta|\cdot C
$$
for some constant $C>0$ which is independent of the choice of $\epsilon$ but may depend on
the choice of the neighborhoods $U_p$, the K\"{a}hler form $\Omega$, the cutoff function
$\rho$ as well as the $1$-form $\gamma$ in the construction of $\tilde{\omega}$. 
With this understood, note that by taking $\epsilon >0$ sufficiently small, we can make
the above difference arbitrarily close to zero, once we fix the other various choices.
On the other hand,  by choosing the neighborhood $U$ in Theorem 1.1 sufficiently small, 
we can also arrange so that the difference
$|c_1(K_{X})\cdot [\omega^\prime]-c_1(K_{X})\cdot [\omega]|$ is arbitrarily close to zero.
With this understood, the difference 
$$|c_1(K_{M_G})\cdot [\tilde{\omega}]-(c_1(K_X)\cdot [\omega]+
\sum_i \frac{1-m_i}{m_i} \omega^\prime(\Sigma_i)+\sum_j\sum_{k\in I_j} a_{j,k} \tilde{\omega}(E_{j,k}))|
$$
can be made arbitrarily close to zero by choosing $\tilde{\omega}$ properly. We claim that
under the assumptions of the lemma, $c_1(K_{M_G})\cdot [\tilde{\omega}]<0$ for such an 
$\tilde{\omega}$. To see this, consider first the case where $c_1(K_M)\cdot [\omega]<0$.
In this case, $c_1(K_X)\cdot [\omega]=\frac{1}{|G|} c_1(K_M)\cdot [\omega]<0$. With the 
other two terms $\sum_i \frac{1-m_i}{m_i} \omega^\prime(\Sigma_i)$ and 
$\sum_j\sum_{k\in I_j} a_{j,k} \tilde{\omega}(E_{j,k}))$ being non-positive, it follows immediately 
that $c_1(K_{M_G})\cdot [\tilde{\omega}]<0$. On the other hand, when $M$ has torsion canonical
class, $c_1(K_X)\cdot [\omega]=\frac{1}{|G|} c_1(K_M)\cdot [\omega]=0$. If the singular set of $M/G$ contains a $2$-dimensional component, then the term 
$\sum_i \frac{1-m_i}{m_i} \omega^\prime(\Sigma_i)$ is negative, which implies that
$c_1(K_{M_G})\cdot [\tilde{\omega}]<0$. If the singular set of $M/G$ contains no 
$2$-dimensional components, then by Proposition 3.2, 
$$
c_1(K_{M_G})\cdot [\tilde{\omega}]=\sum_j\sum_{k\in I_j} a_{j,k} \tilde{\omega}(E_{j,k}),
$$
which is also negative as in this case, the singular set of $M/G$ must contain an isolated 
non-Du Val singularity, so that one of the coefficients $a_{j,k}$ is negative. This finishes 
off the proof of the lemma.

\end{proof}

\noindent{\bf Proof of Theorem 1.9:}

\vspace{2mm}

By Proposition 3.3, we may assume without loss of generality that the symplectic $G$-manifold
$M$ is minimal. We consider first the case where $\kappa^s(M)=0$. Then by Theorem 1.0 in
\cite{C2}, the above assumption means that $M$ is minimal as a smooth $4$-manifold. With this
understood, the assumption $\kappa^s(M)=0$ is equivalent to $M$ having a torsion canonical 
class (cf. \cite{Li}). 

If the $G$-action is free, then $M_G=M/G$, which also has torsion canonical class. Hence
in this case, $\kappa^s(M_G)=0$. If $M/G$ has only isolated singular points which are all 
Du Val singularities, then $c_1(K_{M_G})$ is torsion by Proposition 3.2, and we have 
$\kappa^s(M_G)=0$ as well. In the remaining case, $M_G$ is rational or ruled by Lemma 4.1.

We claim that when $M_G$ is irrational ruled, $b_1(M_G)=2$ so that $M_G$ is a ruled 
surface over $T^2$. Suppose to the contrary that $b_1(M_G)>2$. Then 
$b_1(M)\geq b_1(M/G)=b_1(M_G)\geq 4$. On the other hand, it is known that $b_1(M)\leq 4$
(cf. \cite{Li,Li1, Bauer}), so that $b_1(M)=b_1(M/G)=4$. In particular, the induced $G$-action on $H^1(M,\Q)$ is trivial. Now we note that $b_1(M)=4$ implies that $M$ is a $\Q$-homology $T^4$
(cf. \cite{Li}), and furthermore, the work of Ruberman-Strle in \cite{Sr} implies that $M$ has the 
same $\Q$-cohomology ring of $T^4$. This gives a contradiction, because the
triviality of the induced $G$-action on $H^1(M,\Q)$ implies that the action is also trivial on 
$H^2(M,\Q)$, and consequently, $b_2^{+}(M/G)=b_2^{+}(M)=3$. But this violates 
$b_2^{+}(M/G)=b_2^{+}(M_G)=1$. Hence the claim. 

It remains to consider the case where $M$ is a rational surface. Denote by $\omega_0$ 
the symplectic structure on $M$. Suppose $M$ is $\C\P^2$ or a Hirzebruch surface. Then
one has in this case 
$c_1(K_M)\cdot [\omega_0]<0$, so that by Lemma 4.1, $M_G$ must be rational (note that
$b_1(M)=0$, so $b_1(M_G)=0$ as well). For the remaining case, i.e., 
$M=\C\P^2\# N\overline{\C\P^2}$ for $N\geq 2$, we need to recall some relevant results 
from \cite{CLW} first. 

There are two possibilities for $M$: either $M$ is monotone, i.e., the class $[\omega_0]$ is
a multiple of $c_1(K_{\omega_0})$, or $M$ is a symplectic $G$-conic bundle. In the former
case, it follows easily that $c_1(K_{\omega_0})\cdot [\omega_0]<0$, so that by Lemma 4.1,
$M_G$ is rational. In the latter case, we need further information about the equivariant 
symplectic cone of $M$. Roughly speaking, let $F$ be the fiber class of the symplectic 
$G$-conic bundle on $M$. Then it is shown in \cite{CLW} that for any sufficiently large
$\delta>0$, the class $-c_1(K_{\omega_0})+\delta F$ can be realized by a $G$-invariant
symplectic structure $\omega$ on $M$, where $c_1(K_\omega)=c_1(K_{\omega_0})$. 
With this understood, observing that $c_1(K_{\omega_0})\cdot F=-2$, we have
$$
c_1(K_\omega)\cdot [\omega]=-c_1(K_{\omega_0})^2-2\delta,
$$
which is negative for sufficiently large $\delta>0$. Hence by Lemma 4.1, $M_G$ is rational 
(note that the diffeomorphism type of $M_G$ depends on the smooth orbifold $M/G$ alone, 
not on the symplectic structure on it). This finishes the proof of Theorem 1.9.


\subsection{The case of $\kappa^s(M)=1$}
Let $(M,\omega)$ be a symplectic $4$-manifold with $\kappa^s=1$, equipped with a finite 
symplectic $G$-action such that $b_2^{+}(M/G)>1$. To simplify the situation, we note that
by Proposition 3.3, we may assume the $G$-manifold $M$ is minimal as far as 
$\kappa^s(M_G)$ is concerned. By Theorem 1.0
in \cite{C2}, this condition is equivalent to the smooth $4$-manifold $M$ being minimal. 
Then note that with this assumption, the condition $\kappa^s(M)=1$ is equivalent to 
$c_1(K_M)\cdot [\omega]>0$ and $c_1(K_M)^2=0$. Finally, for simplicity we assume that
$G=\Z_p$ is cyclic of prime order $p$. 

With the preceding understood, for any given $G$-invariant $\omega$-compatible
almost complex structure $J$ on $M$, there is a finite set of $J$-holomorphic curves $\{C_k\}$ with
multiplicities $n_k>0$, which has the following significance (see \cite{CK1}, Theorem 3.2):
\begin{itemize}
\item $c_1(K_M)=\sum_k n_k C_k$, where the set $\cup_k C_k$ (as well as $\{n_k\}$) is $G$-invariant.
\item Any $2$-dimensional fixed component of $G$ is contained in $\cup_k C_k$, and if a fixed point
$m$ of $G$ is not contained in $\cup_k C_k$, then the $G$-action on the tangent space $T_m M$ is
contained in $SL_2(\C)$, in particular, $m$ is an isolated fixed point.
\item If $C_k$ is a fixed component of $G$, then $n_k\geq p-1$ (cf. \cite{C1}, Lemma 1.6). Moreover,
$C_k$ is either a torus of self-intersection zero or a $(-2)$-sphere. In the latter case, $C_k$
must intersect another $(-2)$-sphere in $\cup_k C_k$ which is not fixed by $G$.
\end{itemize}

\vspace{2mm}

\noindent{\bf Proof of Theorem 1.10:}

\vspace{2mm}

We set $X:=M/G$ to be the quotient orbifold. 

For (1), we pick a $G$-invariant $\omega$-compatible almost complex structure $J$, so that
$c_1(K_M)=\sum_k n_k C_k$ for a finite set of $J$-holomorphic curves $\{C_k\}$ with multiplicities $n_k>0$. Then since any $2$-dimensional fixed component of $G$ is contained in $\cup_k C_k$, and moreover, if $C_k$ is a fixed component of $G$, then $n_k\geq p-1$, it follows easily from Proposition 3.2 that $c_1(K_{|X|})$ is torsion if and only if every $C_k$ is fixed by $G$ and $n_k=p-1$ for every $k$. Note that in particular, each $C_k$
must be a torus of self-intersection zero. Furthermore, if every $C_k$ is fixed by $G$, then any
isolated fixed point of $G$ must be in the complement of $\cup_k C_k$, so that the induced 
$G$-action on the tangent space of the fixed point is contained in $SL_2(\C)$. This in particular
implies that $c_1(K_{M_G})$ is torsion. Part (1) follows easily.

For (2), we first note that $c_1(K_{M_G})^2=c_1(K_{|X|})^2+\sum_m K_m^2$, where the second 
term is the sum over all isolated fixed points $m$ of $G$. Thus it suffices to show that 
$c_1(K_{|X|})^2=-\frac{2(p-1)^2}{p}\cdot s$, where 
$s$ is the number of $(-2)$-spheres fixed by $G$. To this end, we denote by $Y$ any 
$2$-dimensional fixed component of $G$, and notice that by Proposition 3.2, 
$c_1(K_X)=c_1(K_{|X|})+\sum_Y \frac{p-1}{p} Y$, where $Y$ is regarded as a surface in $|X|$. Then since $c_1(K_X)^2=\frac{1}{p}c_1(K_M)^2=0$, we have
$$
c_1(K_{|X|})^2+\sum_Y \frac{(p-1)^2}{p} Y^2+\sum_Y  \frac{2(p-1)}{p} c_1(K_{|X|})\cdot Y=0,
$$
where $Y^2$ denotes the self-intersection of $Y$ as a surface in $M$. Here we used 
the fact that as a surface in $|X|$, the self-intersection of $Y$ equals $p\cdot Y^2$. Now by
the adjunction formula, $c_1(K_{|X|})\cdot Y+p\cdot Y^2=2g_Y-2$, where $g_Y$ denotes 
the genus of $Y$. It follows easily that
$$
c_1(K_{|X|})^2=\sum_Y  \frac{2(p-1)}{p}(2-2g_Y)+\sum_Y \frac{p^2-1}{p} Y^2.
$$
Finally, recall that as a surface in $M$, $Y$ is either a torus of self-intersection zero or a
$(-2)$-sphere, which gives immediately that $c_1(K_{|X|})^2=-\frac{2(p-1)^2}{p}\cdot s$.

For (3), we first note that 
$$
\chi(M_G)=\chi(M/G)+\sum_m \chi_m 
$$
and 
$$
\sigma(M_G)=\sigma(M/G)-\sum_m \chi_m,
$$
where $\chi_m$ is the number of exceptional divisors in the minimal resolution of the singular point of $M/G$ corresponding to the isolated fixed point $m$ of $G$.
Combining these two equations, we obtain 
$$
c_1(K_{M_G})^2=2\chi(M_G)+3\sigma(M_G)=2\chi(M/G)+3\sigma(M/G)-\sum_m \chi_m.
$$
On the other hand, when the $G$-action is homologically trivial, $\chi(M/G)=\chi(M)$
and $\sigma(M/G)=\sigma(M)$, so that 
$$
2\chi(M/G)+3\sigma(M/G)=2\chi(M)+3\sigma(M)=c_1(K_M)^2=0.
$$
It follows immediately that $c_1(K_{M_G})^2=-\sum_m \chi_m$.

For (4), we fix a symplectic structure $\tilde{\omega}$ on $M_G$ as constructed in Theorem 1.5.
With this understood, first note that if $c_1(K_{M_G})^2=0$, then the only $2$-dimensional fixed components are tori of self-intersection zero, and for all the isolated fixed points $m$, $K_m=0$. 
This gives immediately $c_1(K_{M_G})=\pi^\ast c_1(K_{|X|})$ by Proposition 3.2. On the other 
hand, if $M_G$ is not minimal, then for any $\tilde{\omega}$-compatible almost complex
structure $J$, there exists a $J$-holomorphic $(-1)$-sphere $C$ in $M_G$ (cf. \cite{C1}, 
Lemma 2.3). We will get a contradiction by showing $\pi^\ast c_1(K_{|X|})\cdot C\geq 0$, because 
$c_1(K_{M_G})\cdot C=-1$. 

To this end, we shall particularly choose an $\tilde{\omega}$-compatible almost complex
structure on $M_G$ as follows. We begin by fixing a $G$-invariant $J$ on $M$ with the
following property: near each of the isolated fixed points $m$, we may identify $\omega, J$ with
the standard symplectic structure and complex structure on $\C^2$, such 
that the $G$-action is given
by a complex linear action (this is possible by the equivariant Darboux theorem). In particular, 
$J$ is integrable near each $m$. With this understood, we have $c_1(K_M)=\sum_k n_k C_k$ 
for a finite set of $J$-holomorphic curves $\{C_k\}$ as
shown in \cite{CK1}. Since $\cup_k C_k$ is $G$-invariant, its image in $X=M/G$ under the 
$G$-action is also a set of $J$-holomorphic curves, which will be denoted by $\{C_i^\prime\}$, with 
multiplicities $n_i^\prime$. We remark that if $C_i^\prime$ is the image of some $C_k$ which is not
fixed under $G$, then $n_i^\prime=n_k$, and if $C_k$ is fixed under $G$, then $n_i^\prime=
n_k/p$. With this understood, note that $c_1(K_X)=\sum_i n_i^\prime C_i^\prime$. Furthermore, since $c_1(K_{M_G})^2=0$, $C_i^\prime$ must be a torus with self-intersection zero if it is 
the image of some $C_k$ which is fixed by $G$, and such a $C_i^\prime$ is disjoint from the other components in $\cup_i C_i^\prime$. In particular, this allows us to modify $J$ near each of such $C_i^\prime$
and obtain an $\omega^\prime$-compatible $J^\prime$ on $|X|$, such that each $C_i^\prime$ in
$\cup_i C_i^\prime$ is $J^\prime$-holomorphic. With this understood, note that by Proposition 3.2, 
$$
c_1(K_{|X|})=\sum_i \hat{n_i} C_i^\prime,
$$
where $\hat{n_i}=n_i^\prime$ if $C_i^\prime$ is not the image of a fixed component, and 
$\hat{n_i}=n_i^\prime-(p-1)/p$ otherwise. Note that $\hat{n_i}\geq 0$ for all $i$.

Now with $|X|$ equipped with $\omega^\prime$ and $J^\prime$, we proceed to the 
construction of $M_G$ with a symplectic structure $\tilde{\omega}$ as in Theorem 1.5.
The key observation here is that since $J^\prime=J$ is integrable near each singular point
of $|X|$, it can be naturally extended to an almost complex structure, still denoted by $J^\prime$,
on $M_G$ which is integrable near the exceptional divisors. Furthermore, $J^\prime$ is
also $\tilde{\omega}$-compatible except in the interpolation regions $A$, where it is only
$\tilde{\omega}$-tame. With this understood, we modify $J^\prime$ in the regions $A$ to 
obtain an $\tilde{\omega}$-compatible $J$ on $M_G$, such that (the proper transform of)
each $C_i^\prime$ is $J$-holomorphic (this is possible because $(\cup_i C_i^\prime)\cap A$ is a disjoint union of embedded symplectic surfaces). This $J$ is the particular $\tilde{\omega}$-compatible almost complex structure on $M_G$ we choose to work with.

With the preceding understood, we now consider the $J$-holomorphic $(-1)$-sphere $C$
in $M_G$. Since $J$ is integrable near  the exceptional divisors and the map $\pi: M_G
\rightarrow |X|$ is simply holomorphically contracting the exceptional divisors to the singular
points in $|X|$, it follows easily that $J$ descends to $|X|$ such that each $C_i^\prime$ in
$|X|$ is $J$-holomorphic. Moreover, the image of $C$ under the map $\pi$, which is 
denoted by $C^\prime$, is also $J$-holomorphic. With this understood, we note that 
$$
\pi^\ast c_1(K_{|X|})\cdot C=c_1(K_{|X|})\cdot C^\prime=(\sum_i \hat{n_i}C_i^\prime)\cdot C^\prime.
$$

To see that $\pi^\ast c_1(K_{|X|})\cdot C\geq 0$, we first consider the case where  for any $i$,
$C^\prime\neq C_i^\prime$. In this case, $\pi^\ast c_1(K_{|X|})\cdot C=(\sum_i \hat{n_i}C_i^\prime)\cdot C^\prime \geq 0$ by the positivity of  intersection of $J$-holomorphic curves. 
If $C^\prime=C_i^\prime$ for some $i$, we will need to recall some additional information
from \cite{CK1} about the curves $\{C_k\}$ in $c_1(K_M)=\sum_k n_k C_k$. It is shown in
Lemma 3.3 of \cite{CK1} that $c_1(K_M)\cdot C_k=0$ for each $k$. From this it follows 
easily that $c_1(K_X)\cdot C_i^\prime=0$ for each $i$, where 
$c_1(K_X)=\sum_i n_i^\prime C_i^\prime$. Now recall that 
$c_1(K_{|X|})=\sum_i \hat{n_i}C_i^\prime$ where $\hat{n_i}=n_i^\prime$ if $C_i^\prime$ is not the image of a fixed component, and if $C_i^\prime$ is the image of a fixed component, we have
$\hat{n_i}=n_i^\prime-(p-1)/p$. Since in the latter case, $C_i^\prime$ is a torus of self-intersection zero which is disjoint from the other components of $\{C_i^\prime\}$, it follows easily that 
$c_1(K_{|X|})\cdot C_i^\prime=0$ for each $i$ as well. This shows that if $C^\prime=C_i^\prime$ for some $i$, then $\pi^\ast c_1(K_{|X|})\cdot C=c_1(K_{|X|})\cdot C^\prime=0$. This finishes the
proof of (4), and the proof of Theorem 1.10 is complete.

\begin{rem}
We shall point out that the expression 
$$
c_1(K_{M_G})^2=-\frac{2(p-1)^2}{p}\cdot s+\sum_m K_m^2
$$
in Theorem 1.10(2) can be also derived independently from the Lefschetz fixed-point theorem 
and the
$G$-signature theorem; in particular, the integrability of the right-hand side does not give any
new constraints to the fixed-point set structure. To see this, note that we have already seen that
$$
c_1(K_{M_G})^2=2\chi(M/G)+3\sigma(M/G)-\sum_m \chi_m.
$$
Now by the Lefschetz fixed-point theorem, we have 
$$
p\cdot \chi(M/G)=\chi(M)+(p-1)\cdot (\sum_m 1+\sum_Y (2-2g_Y)),
$$
and by the $G$-signature theorem, we have
$$
p\cdot \sigma(M/G)=\sigma(M)+\sum_m def_m +\sum_Y def_Y,
$$
where $def_m$ and $def_Y$ stand for the signature defect at $m$ and $Y$ respectively 
(cf. \cite{HZ}).
With this understood, it follows easily that
$$
c_1(K_{M_G})^2=-\sum_m \chi_m+\sum_m \frac{2(p-1)}{p}+\sum_Y \frac{2(p-1)}{p}(2-2g_Y)
+\frac{3}{p}(\sum_m def_m +\sum_Y def_Y).
$$
On the other hand, it is known that if an isolated fixed point $m$ is of type $(1,q)$ for some $0<q<p$,
then $K_m^2+\chi_m= \frac{2(p-1)}{p}-12\cdot s(q,p)$, where $s(q,p)$ denotes the corresponding
Dedekind sum (see e.g. \cite{NN}, \S 7.1 in page 304). Now with the fact that 
$def_m=-4p\cdot s(q,p)$ and $def_Y=\frac{p^2-1}{3}Y^2$ (cf. \cite{HZ}), we see easily that
$$
c_1(K_{M_G})^2=\sum_m K_m^2+ \sum_Y  \frac{2(p-1)}{p}(2-2g_Y)+\sum_Y \frac{p^2-1}{p} Y^2,
$$
which implies easily the expression $c_1(K_{M_G})^2=-\frac{2(p-1)^2}{p}\cdot s+\sum_m K_m^2$.
\end{rem}

\begin{exm}
In this example, we shall examine homologically trivial symplectic $G$-actions on a
symplectic homotopy $K3$ surface $M$ (we assume $\kappa^s(M)=1$; otherwise there are
no such actions, cf. \cite{CK1}). This is one of the main rigidity problems concerning
symplectic finite group actions, so it is interesting to test the strength of our new approach on this
question. For simplicity we assume $G=\Z_3$ (for $G=\Z_2$, there are no homologically 
trivial actions, cf. \cite{Ma, Ru}). 

First of all, we recall that the $2$-dimensional fixed components of $G$
are either tori with self-intersection zero
or $(-2)$-spheres. The former type of components make no contributions in any of the 
$G$-index theorems' calculation; in particular, they can not be detected by these theorems.
For simplicity, we shall ignore the $2$-dimensional toroidal fixed components in the discussion.

With this understood, let $x$ be the number of $(-2)$-spheres fixed by $G$, and let $y,z$ be the
number of isolated fixed points of type $(1,1)$ and $(1,2)$ respectively. Notice that if $m$ is
an isolated fixed point of type $(1,1)$, one has $K_m^2=-1/3$ and $\chi_m=1$, and if $m$ is
of type $(1,2)$, one has $K_m^2=0$ and $\chi_m=2$. Now since the $G$-action is homologically
trivial, the Lefschetz fixed-point theorem gives $\chi(M^G)=\chi(M)=24$. On the other hand,
by Theorem 1.10(2) and (3), we have 
$$
-\frac{2(p-1)^2}{p}\cdot s+\sum_m K_m^2=-\sum_m \chi_m.
$$
It follows easily that $x,y,z$ satisfy the following equations
$$
2x+y+z=24, \;\; \frac{8}{3}x+\frac{1}{3}y=y+2z.
$$
There are four possible solutions, which are listed below:
\begin{itemize}
\item [{(1)}] $x=4$, $y=16$, $z=0$,
\item [{(2)}] $x=5$, $y=11$, $z=3$,
\item [{(3)}] $x=6$, $y=6$, $z=6$,
\item [{(4)}] $x=7$, $y=1$, $z=9$.
\end{itemize}
Case (4) can be further eliminated using the constraints from \cite{CK1}. More concretely,
recall that $c_1(K_M)$ is represented by $J$-holomorphic curves $\{C_k\}$ with multiplicities 
$n_k$. In order to eliminate (4), we appeal to the following properties of the curves $\{C_k\}$:
the connected components of $\cup_k C_k$ fall into five different types, and when $p=3$, one 
can easily see that whenever a connected component of $\cup_k C_k$ contains a 
$(-2)$-sphere fixed by $G$, it also contains an isolated fixed point of type $(1,1)$ (cf. \cite{CK1},
Proposition 3.7). In particular, this implies $y\geq x$. The remaining three cases can not
be eliminated by the constraints from \cite{CK1}, neither do they violate any known obstructions 
for smoothable $\Z_p$-actions (e.g. as listed in \S 3 of \cite{CK2}).

With the preceding understood, we look at the resolution $M_G$ of these remaining actions
listed in (1)-(3).
From the proof of Theorem 1.10, we have seen how to determine $c_1(K_{M_G})$ from a
concrete $J$-holomorphic representative of $c_1(K_M)$, i.e, $c_1(K_M)=\sum_k n_k C_k$
from \cite{CK1}. For illustration we shall examine one particular case in details; all other 
possibilities (there are finitely many of them) can be similarly analyzed and they all give the
same conclusion. 

Let's assume $c_1(K_M)=\sum_k n_k C_k$, where the connected components of $\cup_k C_k$
consist of four tori of self-intersection zero, denoted by $T_i$, $i=1,2,3,4$, with multiplicities 
$\alpha_i$, and four other connected components, denoted by $\Lambda_j$, $j=1,2,3,4$, where
each $\Lambda_j$ is a union of three $(-2)$-spheres whose intersection graph forms a cycle 
(it is denoted by $\tilde{A}_2$ in \cite{CK1}). We denote by $\beta_j$ the multiplicity of each
$(-2)$-sphere in $\Lambda_j$. Concerning the fixed-point set, each $T_i$ contains three 
isolated fixed points of type $(1,1)$, and each $\Lambda_j$ contains a fixed $(-2)$-sphere and
an isolated fixed point of type $(1,1)$, which is the intersection point of the two $(-2)$-spheres 
in $\Lambda_j$ that are not fixed by $G$ . So totally, we have sixteen isolated 
fixed points of type $(1,1)$, four fixed $(-2)$-spheres, and no isolated fixed points of type $(1,2)$
(i.e., we are in case (1), where $x=4$, $y=16$, $z=0$). 

With this understood, we observe that each $T_i$ is mapped to an orbifold two-sphere in $M/G$ 
which becomes a $(-1)$-sphere in $M_G$ (i.e., the proper transform) after resolving the singularities
of $M/G$. We denote this $(-1)$-sphere by $S_i$. It follows easily that each $T_i$ gives rise to
a configuration of four two-spheres in $M_G$: one  $(-1)$-sphere, $S_i$, and three $(-3)$-spheres,
denoted by $E_{i,1}, E_{i,2},E_{i,3}$, which are the corresponding exceptional divisors in the
resolution. On the other hand, each $\Lambda_j$ also gives rise to a configuration of four 
two-spheres in $M_G$: one $(-6)$-sphere, denoted by $A_j$ which is the fixed $(-2)$-sphere
in $\Lambda_j$, two $(-1)$-spheres, denoted by $B_{j,1}$, $B_{j,2}$, which are the image of
the two $(-2)$-spheres in $\Lambda_j$ that are not fixed by $G$, and one $(-3)$-sphere,
denoted by $E_j$, which is the exceptional sphere from the isolated fixed point in $\Lambda_j$.
With this understood, we note that
$$
c_1(K_{M_G})=\sum_{i=1}^4(\alpha_i S_i+ \frac{\alpha_i-1}{3}\sum_{k=1}^3 E_{i,k})+
\sum_{j=1}^4(\frac{\beta_j-2}{3}A_j+\beta_j(B_{j,1}+B_{j,2})+\frac{2\beta_j-1}{3}E_j).
$$

We see immediately twelve disjoint $(-1)$-spheres in $M_G$, which are $S_i$, $i=1,2,3,4$,
and $B_{j,1},B_{j,2}$, $j=1,2,3,4$. Blow down these $(-1)$-spheres and denote the images
of $E_{i,k}$, $E_j$  by $E_{i,k}^\prime$, $E_j^\prime$ respectively. Then it is easy to see that
each $E_{i,k}^\prime$ is a $(-2)$-sphere and each $E_j^\prime$ is a $(-1)$-sphere. 
Further blow down each $E_j^\prime$, we arrive at a symplectic $4$-manifold, denoted by
$M_G^\prime$. We observe that each $(-6)$-sphere $A_j$ becomes a nodal $(-2)$-sphere in
$M_G^\prime$, which we denote by $A_j^\prime$. With this understood, we note that
$$
c_1(K_{M_G^\prime})=\sum_{i=1}^4\frac{\alpha_i-1}{3} (E_{i,1}^\prime+E_{i,2}^\prime
+E_{i,3}^\prime)+\sum_{j=1}^4\frac{\beta_j-2}{3}A_j^\prime.
$$
Since each $E_{i,k}^\prime$, $A_j^\prime$ is $J$-holomorphic for some $J$ on $M_G^\prime$,
and none of them is a $(-1)$-sphere, it follows that $M_G^\prime$ is minimal (see \cite{C1}, 
Lemma 2.3). On the other hand, $c_1(K_{M_G})^2=-16$, and we have blown down 
successively sixteen $(-1)$-spheres in $M_G$ to reach $M_G^\prime$, it follows easily
that $c_1(K_{M_G^\prime})^2=0$. This shows that $\kappa^s(M_G)\leq 1$. 

In conclusion, one can verify that Conjecture 1.7 is true for homologically trivial symplectic 
$\Z_3$-actions on a symplectic homotopy $K3$ surface. On the other hand, for these group 
actions our new approach does not give any new constraints to the fixed-point set structure.

\end{exm}

{\bf Acknowledgements:} 
We are extremely grateful to R. Inanc Baykur for his interest in this work, and particularly, for 
pointing out a serious error in the previous version of this paper 
(i.e., the preprint arXiv:1708.09428v1 [math.SG]).
The error occurred in the section on symplectic Calabi-Yau surfaces, which is removed from the current version; 
that topic will be addressed elsewhere \cite{C3}. 
Part of the work was done during a visit at Capital Normal University in Beijing. We are 
very grateful to Professor Fuquan Fang for the invitation and the very warm hospitality during
the visit. Finally, we wish to thank the referee for a very thorough reading of the manuscript.

\end{document}